\documentclass{amsart}[12pt]
\usepackage{amssymb,amsmath}
\usepackage{enumerate}
\usepackage[english]{babel}
\usepackage{color}
\usepackage{comment}

\newtheorem{teo}{Theorem}[section]
\newtheorem{pro}{Proposition}[section]
\newtheorem{cor}{Corollary}[section]

\newtheorem{lm}{Lemma}[section]

\theoremstyle{definition}
\newtheorem{rem}{Remark}[section]

\newtheorem{ex}{Example}[section]

\title[Noncommutative Breiman ergodic theorem]{Noncommutative extensions of an ergodic theorem of Breiman}

\author{Morgan O'Brien}
\address{Minnesota State University, Mankato\\ Department of Mathematics and Statistics \\ Mankato, MN 56001, USA}
\email{morgan.obrien.2@mnsu.edu, obrienmorganc@gmail.com}

\begin{document}
\begin{abstract}
Let $(\mathcal{M},\tau)$ be a tracial von Neumann algebra, $T:L_1(\mathcal{M},\tau)\to L_1(\mathcal{M},\tau)$ be induced by a normal $\tau$-preserving $*$-automorphism of $\mathcal{M}$, and $x\in L_1(\mathcal{M},\tau)$. We will prove the averages $$\frac{1}{n}\sum_{k=0}^{n-1}T^k(x_k)$$ converge to $\tau(x)\textbf{1}$ bilaterally almost uniformly as $n\to\infty$, generalizing an ergodic theorem proven by Breiman in \cite{Breiman57}, where $(x_k)_{k=0}^{\infty}\in L_p(\mathcal{M},\tau;\ell_\infty)$ is such that $(x_k-x)_{k=0}^{\infty}\in L_p(\mathcal{M},\tau;c_0)$. A generalization to semifinite $(\mathcal{M},\tau)$ and general positive Dunford-Schwartz operators $T\in DS^+(\mathcal{M},\tau)$ is shown as well when $1<p<\infty$. In addition to these, we also prove a number of generalizations of the results to weighted, subsequential, and other types of averages.
\end{abstract}
\date{November 16, 2026}
\maketitle

\section{Introduction}\label{Intro}

The ergodic theorems of Birkhoff and von Neumann study the various types of convergence of the averages
$$\frac{1}{n}\sum_{k=0}^{n-1}f\circ\phi^k,$$ where $(\Omega,\mathcal{F},\nu)$ is a measure space, $\phi:\Omega\to\Omega$ is a measure-preserving transformation, $1\leq p<\infty$, and $f\in L_p(\Omega,\nu)$, where Birkhoff's individual ergodic theorem shows that these averages converge $\nu$-a.e., while von Neumann's mean ergodic theorem proves the $L_2$-norm convergence of these averages when $p=2$; note that $L_p$-norm convergence does hold when $1\leq p<\infty$.

A major topic in ergodic theory studies in what ways these results can be extended. For example, the terms $f\circ\phi^k$ being averaged over are often replaced by functions of the form $T^k(f)$, where $T:L_p(\Omega,\mathcal{F},\nu)\to L_p(\Omega,\mathcal{F},\nu)$ is a bounded linear operator. Weighted and subsequential ergodic theorems consider the norm and $\nu$-a.e. convergence of averages of the form
$$\frac{1}{n}\sum_{k=0}^{n-1}\alpha_kf\circ\phi^k\text{ and }\frac{1}{n}\sum_{j=0}^{n-1}f\circ\phi^{k_j},$$ 
where $\alpha=(\alpha_k)_{k=0}^{\infty}\subset\mathbb{C}$ and $\textbf{k}=(k_j)_{j=0}^{\infty}\subset\mathbb{N}_0$ is a strictly increasing sequence. Examples of sequences where the convergence of such averages are known include when $\alpha=(\lambda^k)_{k=0}^{\infty}$ for $\lambda\in\mathbb{T}:=\{\lambda\in\mathbb{C}:|\lambda|=1\}$, when $\alpha=(g(\psi^k(\theta)))_{k=0}^{\infty}$ if $(\Theta,\mathcal{G},\rho)$ is a probability space, $\psi:\Theta\to\Theta$ is a measure-preserving transformation, $g\in L_\infty(\Theta,\mathcal{G},\rho)$, and $\theta\in\Theta$, when $\alpha_k=\mu(k)$ with $\mu:\mathbb{N}_0\to\mathbb{C}$ being the M\"{o}bius function, and when $\textbf{k}$ is the set of primes in $\mathbb{N}_0$. Among many other sources containing results of this nature, see \cite{eakpldlr,bfko89,comli2,ek,lot,rn,wie88}. See Section \ref{SupplementaryResults} below as well.

Another interesting variation of these results replace the usage of a single function $f\in L_p(\Omega,\nu)$ with a sequence of functions $(f_k)_{k=0}^{\infty}\subset L_1(\Omega,\nu)$, studying the different types of convergence of averages of the form
\begin{align}\label{Averages2}
\frac{1}{n}\sum_{k=0}^{n-1}T^k(f_k).
\end{align}
In classical ergodic theory a few results regarding these averages are known. 

One frequently studied type of sequence $(f_k)_{k=0}^{\infty}$ is a strongly $p$-bounded $T$-admissible process. Let $\phi:\Omega\to \Omega$ be a measure-preserving transformation on $(\Omega,\mathcal{F},\nu)$ and define $T_\phi(f):=f\circ\phi$ for $f\in L_p(\Omega,\mathcal{F},\nu)$ and $1\leq p<\infty$. Then $(f_k)_{k=0}^{\infty}$ is said to be a $T_\phi$-admissible process if each $f_k$ is real-valued and if $T_\phi(f)\leq f_{k+1}$ for every $k\geq0$, while the sequence is strongly $p$-bounded if $\sup_{k\in\mathbb{N}_0}\|f_k\|_p<\infty$. Then the $\nu$-a.e. convergence of the averages (\ref{Averages2}) is proven in \cite{comli2} (among other sources), even for weighted averages of that form; in fact, a Wiener-Wintner form of the result is shown in that article, allowing the averages to be weighted by $(\lambda^k)_{k=0}^{\infty}$ for each $\lambda\in\mathbb{T}$ with the $\nu$-a.e. convergence set obtained being independent of $\lambda$.

Another such result in the commutative setting was proven by Breiman in \cite{Breiman57}, where such a theorem was used as an important tool in studying the entropy of certain stochastic processes. A few modifications to how the original result is stated will need to be made in order to be able adapt it to the von Neumann algebra setting. Due to this we state the original form of Breiman's result here.

\begin{teo}\label{BreimanOriginalResult}(Cf. \cite[Theorem 1]{Breiman57})
Let $(\Omega,\mathcal{F},\nu)$ be a probability space, $T:\Omega\to\Omega$ be an ergodic measure-preserving transformation, and $(g_k)_{k=0}^{\infty}\subset L_1(\Omega,\mathcal{F},\nu)$ and $g\in L_1(\Omega,\mathcal{F},\nu)$ be such that $\int_{\Omega}\sup_{k\in\mathbb{N}_0}|g_k|d\nu<\infty$ and $g_k\to g$ $\nu$-a.e. Then
$$\frac{1}{n}\sum_{k=0}^{n-1}g_k(T^k(\omega))\to\int_{\Omega}g\,d\mu\text{ as }n\to\infty$$ for $\nu$-a.e. $\omega\in\Omega$.
\end{teo}	

In this article we will extend this result to hold for the noncommutative $L_p$-spaces associated to a von Neumann algebra, as well as prove a number of variations and generalizations to it. The relevant convergence will be bilaterally almost uniform convergence (abbreviated as b.a.u.), which is the adaptation of almost uniform convergence to the von Neumann algebra setting (and so implies almost everywhere convergence in the commutative setting). See Section \ref{Preliminaries} for any unexplained notation and definitions below.

In Section \ref{MainResults} we will study a few direct variants of Breiman's result. In Theorem \ref{DirectBreimanErgodicTheorem} the direct noncommutative analogue of Breiman's result above will be shown for tracial von Neumann algebras. Before that is done, we first also prove a general version for positive Dunford-Schwartz operators in Theorem \ref{GeneralBreimanErgodicTheorem} for semifinite von Neumann algebras $(\mathcal{M},\tau)$ that allows one to consider operators in $L_p(\mathcal{M},\tau)$ for some $1<p<\infty$ to be used. While the proofs for two results will be similar, the proof of the latter result a bit easier due to the setting allowing strong type maximal inequalities to be used, while the former setting requires a bit more detailed analysis to prove. Although our statements of these results look a bit different from Theorem \ref{BreimanOriginalResult}, we prove at the end of that section in Lemma \ref{WhyAssumption} that the formulations are equivalent in the commutative setting.

In Section \ref{SupplementaryResults} we will prove a number of modifications and generalizations of these results. Most of the results in this section will be proven in one of two ways. In one method we will construct new sequences in $L_p(\mathcal{M},\tau;c_0)$ from old ones and apply Theorems \ref{GeneralBreimanErgodicTheorem} and \ref{DirectBreimanErgodicTheorem} to the sequences obtained in this way. Typically this will be done by showing that an appropriate operator from $L_p(\mathcal{M},\tau;\ell_\infty)$ to itself is bounded and leaves the subspace $L_p(\mathcal{M},\tau;c_0)$ invariant. For the other method we will show how to rewrite the results in a way so that the proofs can be obtained by repeating the proofs of Theorems \ref{GeneralBreimanErgodicTheorem} and \ref{DirectBreimanErgodicTheorem}.

For the first extension of Breiman's result we will study weighted versions of it, which will be proven by applying the original result to cleverly chosen sequences in Theorem \ref{WeightedDirect} for bounded sequences as weights, generalizing Theorem \ref{DirectBreimanErgodicTheorem}. Afterwards, in Theorem \ref{WeightedGeneralized} we show that unbounded sequences can be used in the corresponding generalizations of Theorem \ref{GeneralBreimanErgodicTheorem}.

After this we will discuss subsequential ergodic theorems, noting that the same methods used to prove Theorem \ref{GeneralBreimanErgodicTheorem} can also extend to averages along powers coming from appropriate subsequences of $\mathbb{N}$ so long as an appropriate maximal inequality holds. For the case of the powers belonging to the sequence $(k^t)_{k=0}^{\infty}$ for some fixed $t\in\mathbb{N}$ we will also show in Theorem \ref{MobiusPolynomialErgodicTheorem} that the weighted subsequential averages $\frac{1}{n}\sum_{k=0}^{n-1}\mu(k)T^{k^t}(x)$ converge b.a.u. as $n\to\infty$ for every $1<p<\infty$ and $x\in L_p(\mathcal{M},\tau)$ when $T:L_p(\mathcal{M},\tau)\to L_p(\mathcal{M},\tau)$ is induced by a normal $\tau$-preserving $*$-automorphism. 

Another application will follow by applying Theorem \ref{GeneralBreimanErgodicTheorem} to certain sequences of the form $(RS^k(x))_{k=0}^{\infty}$ when $1<p<\infty$, $S\in DS^+(\mathcal{M},\tau)$, $R:L_p(\mathcal{M},\tau)\to L_p(\mathcal{M},\tau)$ is a bounded linear operator, and $x\in L_p(\mathcal{M},\tau)$, with the corresponding averages taking the form $$\frac{1}{n}\sum_{k=0}^{n-1}T^kRS^k(x)$$ for any $T\in DS^+(\mathcal{M},\tau)$, which resembles the entangled ergodic theorems considered by Liebscher in \cite{lieb}; see \cite{ekk} for a more detailed account on the history of entangled ergodic theorems in the commutative setting. In the present article, if $1<p<\infty$ and $x\in L_p(\mathcal{M},\tau)$ then the sequence $(S^k(x))_{k=0}^{\infty}$ is known to converge to $F_{S}(x)$ b.a.u. as $k\to\infty$ when the restriction of $S$ to $L_2(\mathcal{M},\tau)$ is positive as a Hilbert-space operator on $L_2(\mathcal{M},\tau)$, from which applying Theorem \ref{GeneralBreimanErgodicTheorem} would imply that $\frac{1}{n}\sum_{k=0}^{n-1}T^kRS^k(x)$ converges to $F_{T}RF_{S}(x)$ b.a.u. as $n\to\infty$ for every $T\in DS^+(\mathcal{M},\tau)$. In Theorem \ref{EntangledAverages2} we will show a further generalization of this still holds when the assumptions on $S$ above only hold for $S^d$ for some fixed $d\in\mathbb{N}$.

Finally, we will relate the results of this article to other known results concerning $T$-admissible sequences as considered in \cite{comli2} (among other sources) in the commutative setting and \cite{comli1,ob5} in the noncommutative setting. While the results of the present article do apply to $T$-admissible sequences, we discuss how the study of $T$-admissible sequences on their own differs from the results in this article.

\section{Preliminaries}\label{Preliminaries}

Throughout this article we will let $\mathbb{N}$ denote the set of natural numbers, and will write $\mathbb{N}_0=\mathbb{N}\cup\{0\}$. We will let $\mathcal{H}$ denote a complex Hilbert space, and we will denote the operator norm of an element $x$ of $\mathcal{B(H)}$ (the space of bounded linear operators on $\mathcal{H}$) by $\|x\|_\infty$. The identity operator in $\mathcal{B(H)}$ will be denoted by $\textbf{1}$, and if $e\in\mathcal{B(H)}$ is a projection then we will write $e^{\perp}:=\textbf{1}-e$.

The pair $(\mathcal{M},\tau)$ is called a semifinite von Neumann algebra if $\mathcal{M}$ is a von Neumann algebra acting on a Hilbert space $\mathcal{H}$ and $\tau$ is a normal semifinite faithful trace on $\mathcal{M}$. The pair $(\mathcal{M},\tau)$ is called a tracial von Neumann algebra if it is semifinite and if $\tau(\textbf{1})=1$, i.e. if the trace $\tau$ is also a state on $\mathcal{M}$. We will write $\mathcal{P(M)}$ for the set of all projections in $\mathcal{M}$.

Throughout the rest of this section we will let $(\mathcal{M},\tau)$ denote a fixed semifinite von Neumann algebra acting on $\mathcal{H}$.

A closed densely defined operator $x:\mathcal{D}_{x}\to\mathcal{H}$ with $\mathcal{D}_{x}\subset\mathcal{H}$ is said to be affiliated with $\mathcal{M}$ if $xu=ux$ for every unitary operator $u\in\mathcal{M}^{\prime}$, where $\mathcal{M}^{\prime}\subset\mathcal{B(H)}$ denotes the commutant of $\mathcal{M}$. If $x$ is affiliated to $\mathcal{M}$, then it is said to be $\tau$-measurable if for every $\epsilon>0$ there exists $e\in\mathcal{P(M)}$ such that $\tau(e^\perp)\leq\epsilon$ and $xe\in\mathcal{M}$. We will let $L_0(\mathcal{M},\tau)$ denote the space of all $\tau$-measurable operators affiliated with $\mathcal{M}$.

For every $\epsilon,\delta>0$ define the set $V(\epsilon,\delta)$ as
$$\{x\in L_0(\mathcal{M},\tau):\text{ there exists }e\in\mathcal{P(M)}\text{ with }\tau(e^\perp)\leq\epsilon\text{ such that }\|xe\|_\infty\leq\delta\}.$$
The family $\{V(\epsilon,\delta):\epsilon,\delta>0\}$ forms a neighborhood of open sets of $0$ in $L_0(\mathcal{M},\tau)$ which gives rise to the topology of convergence in measure. Equipped with this topology and the operations of closed sum, product, and adjoint, $L_0(\mathcal{M},\tau)$ is a completely metrizable topological $*$-algebra. See \cite{ne} for more details.

An operator $x\in L_0(\mathcal{M},\tau)$ is positive, denoted by $x\geq0$, if $\langle x\xi,\xi\rangle\geq0$ for every $\xi\in\mathcal{D}_{x}$. If $S\subset L_0(\mathcal{M},\tau)$, then we will let $S^{+}:=\{x\in S:x\geq0\}$ and $S^{h}:=\{x\in S:x^*=x\}$.

If $x\in L_0(\mathcal{M},\tau)$ then we can consider its polar decomposition as $x=u|x|$, where $u\in\mathcal{M}$ is such that $\|u\|_\infty\leq1$ and $|x|\in L_0(\mathcal{M},\tau)^{+}$. If $x\in L_0(\mathcal{M},\tau)^{+}$ then we can consider its spectral decomposition as $x=\int_{[0,\infty)}\lambda de_{\lambda}$. Given such a decomposition, we will define
$$\tau(x):=\sup_{n\in\mathbb{N}}\tau\left(\int_{[0,n]}\lambda de_{\lambda}\right).$$ 
This extends the trace $\tau$ from $\mathcal{M}^{+}$ to $L_0(\mathcal{M},\tau)^{+}$.

Given $1\leq p<\infty$ and $x\in L_0(\mathcal{M},\tau)$ we will define 
$$\|x\|_p:=\tau(|x|^{p})^{1/p}$$ and the noncommutative $L_p$-space associated to $(\mathcal{M},\tau)$ to be
$$L_p(\mathcal{M},\tau):=\{x\in L_0(\mathcal{M},\tau):\|x\|_p<\infty\}.$$
We will also denote $L_\infty(\mathcal{M},\tau):=\mathcal{M}$. Then $(L_p(\mathcal{M},\tau),\|\cdot\|_p)$ is a Banach space for every $1\leq p\leq\infty$. It is known that convergence with respect to $\|\cdot\|_p$ implies convergence in measure on $L_0(\mathcal{M},\tau)$. 

If $x\in L_0(\mathcal{M},\tau)^{h}$, then $-|x|\leq x\leq |x|$. If $1\leq p\leq\infty$ and $x\in L_0(\mathcal{M},\tau)$ then $\|x\|_p=\||x|\|_p$, where $x=u|x|$ is the polar decomposition.

Fix $1\leq p<\infty$. We will write $L_p(\mathcal{M},\tau;\ell_\infty)$ for the set of all sequences $(x_k)_{k=0}^{\infty}\subset L_p(\mathcal{M},\tau)$ such that there exists $a,b\in L_{2p}(\mathcal{M},\tau)$ and $(y_k)_{k=0}^{\infty}\in\mathcal{M}$ with $\sup_{k\in\mathbb{N}_0}\|y_k\|_\infty<\infty$ such that
$$x_k=ay_kb\text{ for every }k\in\mathbb{N}_0.$$
Equipped with the norm
$$\|(x_k)_{k=0}^{\infty}\|_{L_p(\mathcal{M},\tau;\ell_\infty)}=\inf_{x_k=ay_kb}\left\{\|a\|_{2p}\|b\|_{2p}\sup_{k\in\mathbb{N}_0}\|y_k\|_\infty\right\},$$ where the infimum is taken over all factorizations mentioned above, the space $L_p(\mathcal{M},\tau;\ell_\infty)$ is a Banach space under coordinate-wise operations; see \cite[Theorem 3.2]{dj} for details. Given $(x_k)_{k=0}^{\infty}\in L_p(\mathcal{M},\tau;\ell_\infty)$ we will adopt the frequently used notation 
$$\left\|{\sup_{k\in\mathbb{N}_0}}^{+}x_k\right\|_p:=\|(x_k)_{k=0}^{\infty}\|_{L_p(\mathcal{M},\tau;\ell_\infty)},$$ noting that the expression ${\sup_{k\in\mathbb{N}_0}}^{+}x_k$ is used for intuition and does not necessarily have any meaning on its own.

We will sometimes index an element of $L_p(\mathcal{M},\tau;\ell_\infty)$ by $\mathbb{N}$ instead of $\mathbb{N}_0$. In this case, if we write $(x_n)_{n=1}^{\infty}\in L_p(\mathcal{M},\tau;\ell_\infty)$, then we view this as reindexing via $(x_n)_{n=1}^{\infty}=(x_{k+1})_{k=0}^{\infty}$. This should not cause any confusion given the context of where the different index sets are used, and we will usually use the letter $k$ for elements of $\mathbb{N}_0$ and $n$ for elements in $\mathbb{N}$.

The following is a very useful characterization of self-adjoint elements of the space $L_p(\mathcal{M},\tau;\ell_\infty)$ which will be used frequently throughout this article.

\begin{lm}(Cf. \cite[Lemma 2.1]{lmz})
Let $1\leq p<\infty$ and $(x_k)_{k=0}^{\infty}\subset L_p(\mathcal{M},\tau)$ be such that $x_k^*=x_k$ for every $k\in\mathbb{N}_0$. 
\begin{itemize}
\item[(i)] If $(x_k)_{k=0}^{\infty}\in L_p(\mathcal{M},\tau;\ell_\infty)$ with $\|{\sup_{k\in\mathbb{N}_0}}^{+}x_k\|_p<K$ for some $K>0$, then there exists $c\in L_p(\mathcal{M},\tau)^{+}$ such that
$$\|c\|_p<K\text{ and }-c\leq x_k\leq c\text{ for every }k\in\mathbb{N}_0.$$
\item[(ii)] If there exists $c,d\in L_p(\mathcal{M},\tau)^{+}$ such that $-d\leq x_k\leq c$ for every $k\in\mathbb{N}_0$, then $(x_k)_{k=0}^{\infty}\in L_p(\mathcal{M},\tau;\ell_\infty)$ with $\|{\sup_{k\in\mathbb{N}_0}}^{+}x_k\|_p\leq\|c+2d\|_p$.
\end{itemize}
\end{lm}

We define $L_p(\mathcal{M},\tau;c_0)$ to be the set of all $(x_k)_{k=0}^{\infty}\in L_p(\mathcal{M},\tau;\ell_\infty)$ which has a factorization $x_k=ay_kb$ such that $a,b\in L_{2p}(\mathcal{M},\tau)$ and $(y_k)_{k=0}^{\infty}\subset\mathcal{M}$ is such that $\lim_{k\to\infty}\|y_k\|_\infty=0$. It is known that $L_p(\mathcal{M},\tau;c_0)$ is a closed subspace of $L_p(\mathcal{M},\tau;\ell_\infty)$.

The following characterization of self-adjoint elements of $L_p(\mathcal{M},\tau;c_0)$ will also be used frequently in this article.

\begin{lm}(Cf. \cite[Corollary 2.3]{lmz})
	Let $1\leq p<\infty$ and $(x_k)_{k=0}^{\infty}\subset L_p(\mathcal{M},\tau)$ be such that $x_k^*=x_k$ for every $k\in\mathbb{N}_0$. Then $(x_k)_{k=0}^{\infty}\in L_p(\mathcal{M},\tau;c_0)$ if and only if for every $\epsilon>0$ there exists $k^{\prime}\in\mathbb{N}_0$ and $c,d\in L_p(\mathcal{M},\tau)^{+}$ such that
	$$\max\{\|c\|_p,\|d\|_p\}\leq\epsilon\text{ and }-d\leq x_k\leq c\text{ for every }k\geq k^{\prime}.$$
\end{lm}

We will write $L_p(\mathcal{M},\tau;c_{00})$ for the subset of all finite sequences in $L_p(\mathcal{M},\tau;\ell_\infty)$. In other words, $(x_k)_{k=0}^{\infty}\in L_p(\mathcal{M},\tau;c_{00})$ if there exists $N\in\mathbb{N}_0$ such that $x_k=0$ for $k\geq N$. The main property of this space that we will need is that $L_p(\mathcal{M},\tau;c_{00})$ is dense in $L_p(\mathcal{M},\tau;c_0)$ with respect to the norm $\|\cdot\|_{L_p(\mathcal{M},\tau;\ell_\infty)}$.

Let $(x_k)_{k=0}^{\infty}\subset L_0(\mathcal{M},\tau)$ and $x\in L_0(\mathcal{M},\tau)$. Then we say that $x_k\to x$ bilaterally almost uniformly (abbreviated as b.a.u.) as $k\to\infty$ if for every $\epsilon>0$ there exists $e\in\mathcal{P(M)}$ such that
$$\tau(e^{\perp})\leq\epsilon\text{ and }\|e(x_k-x)e\|_\infty\to0\text{ as }k\to\infty.$$

It was shown in \cite{cls} that if $x_k\to x$ b.a.u. as $k\to\infty$ then $x_k\to x$ with respect to the measure topology on $L_0(\mathcal{M},\tau)$. In the commutative setting b.a.u. convergence corresponds with almost uniform convergence, which implies almost everywhere convergence (on a finite measure space it is actually equivalent to almost everywhere convergence by Egorov's theorem).

The following is a well-known and well-used result which relates b.a.u. convergence and elements of $L_p(\mathcal{M},\tau;c_0)$.

\begin{lm}(Cf. \cite[Lemma 6.2]{jx})
	Let $1\leq p<\infty$ and $(x_k)_{k=0}^{\infty}\subset L_p(\mathcal{M},\tau)$. 
	If $(x_k)_{k=0}^{\infty}\in L_p(\mathcal{M},\tau;c_0)$, then $x_k\to0$ b.a.u. as $k\to\infty$.
\end{lm}

An operator $T:L_1(\mathcal{M},\tau)+\mathcal{M}\to L_1(\mathcal{M},\tau)+\mathcal{M}$ is called a positive Dunford-Schwartz operator if
$$\|T(x)\|_p\leq\|x\|_p\text{ for every }x\in L_p(\mathcal{M},\tau)\text{ and }1\leq p\leq\infty$$ and if $T(x)\geq0$ for every $x\in (L_1(\mathcal{M},\tau)+\mathcal{M})^{+}$. We will denote the set of all such operators by $DS^+(\mathcal{M},\tau)$. We will call such an operator $T$ ergodic if $T(x)=x$ for some $x\in L_2(\mathcal{M},\tau)$ implies that $x=c\textbf{1}$ for some $c\in\mathbb{C}$ (with $c=0$ if $\tau(\textbf{1})=\infty$).

If $1\leq p<\infty$, $T\in DS^+(\mathcal{M},\tau)$, and $n\in\mathbb{N}$ are all fixed, then for every $x\in L_p(\mathcal{M},\tau)$ we will define the $n$-th ergodic average of $x$ under $T$ by
$$M_n(T)(x):=\frac{1}{n}\sum_{k=0}^{n-1}T^k(x).$$
When $1<p<\infty$, for every $T\in DS^+(\mathcal{M},\tau)$ we will write $F_T$ for the projection of $L_p(\mathcal{M},\tau)$ onto the fixed point subspace of $T$ in $L_p(\mathcal{M},\tau)$; note then that $T(F_T(x))=F_T(x)$ for every $x\in L_p(\mathcal{M},\tau)$. Although some reference to $p$ should be made in the notation for $F_T$, if $1<r<\infty$ and $x\in L_p(\mathcal{M},\tau)\cap L_r(\mathcal{M},\tau)$ then the value of $F_T(x)$ when obtained from $F_T$ as a map from $L_p(\mathcal{M},\tau)$ to itself agrees with what would be found from $F_T$ mapping $L_r(\mathcal{M},\tau)$ to itself. As such, no reference to $p$ needs to be made when writing $F_T$.

If $1\leq p<\infty$ is fixed, and if $T_k:L_p(\mathcal{M},\tau)\to L_p(\mathcal{M},\tau)$ is a linear operator for every $k\in\mathbb{N}_0$, then the family $(T_k)_{k\in\mathbb{N}_0}$ is said to be of strong type $(p,p)$ with constant $C$ if the linear operator $\overline{T}:L_p(\mathcal{M},\tau)\to L_p(\mathcal{M},\tau;\ell_\infty)$ defined by
$$\overline{T}(x):=(T_k(x))_{k=0}^{\infty}$$ is bounded with an operator norm of at most $C$.

For convenience we state the following result which will play a major role throughout this article.

\begin{teo}(Cf. \cite[Theorems 4.1 and 6.3]{jx})
Let $(\mathcal{M},\tau)$ be a semifinite von Neumann algebra, $1<p<\infty$, and $T\in DS^+(\mathcal{M},\tau)$. Then $(M_n(T))_{n\in\mathbb{N}}$ is of strong type $(p,p)$, and for every $x\in L_p(\mathcal{M},\tau)$ it follows that
$$\left(M_n(T)(x)-F_T(x)\right)_{n=1}^{\infty}\in L_p(\mathcal{M},\tau;c_0),$$ where $F_T$ denotes the projection of $L_p(\mathcal{M},\tau)$ onto the fixed point subspace of $T$.
\end{teo}

\section{Noncommutative analogues of Breiman's Ergodic Theorem}\label{MainResults}

In this section we will prove a noncommutative version of the ergodic theorem of Breiman considered in \cite{Breiman57}, as well as a generalization of it to general noncommutative $L_p$-spaces when $1<p<\infty$. We prove the latter result first since, while the methods used to prove both results will be similar, the proof for the case $1<p<\infty$ is a bit simpler due to the availability of strong type maximal inequalities in this range.

%\subsection{Standard Ergodic Averages}

Compared to the original statement of Theorem \ref{BreimanOriginalResult}, in all of our results we will replace the assumptions relating to $g_k$ and $g$ in Theorem \ref{BreimanOriginalResult} by instead assuming that $(x_k)_{k=0}^{\infty}\in L_p(\mathcal{M},\tau;\ell_\infty)$ and $x\in L_p(\mathcal{M},\tau)$ are such that $(x_k-x)_{k=0}^{\infty}\in L_p(\mathcal{M},\tau;c_0)$. We justify this at the end of this section in Lemma \ref{WhyAssumption}, where it is shown in the commutative setting that the assumptions of Theorem \ref{BreimanOriginalResult} imply this different looking assumption.

\begin{teo}\label{GeneralBreimanErgodicTheorem}
Let $(\mathcal{M},\tau)$ be a semifinite von Neumann algebra, $T\in DS^+(\mathcal{M},\tau)$, and $1<p<\infty$. If $(x_k)_{k=0}^{\infty}\in L_p(\mathcal{M},\tau;c_0)$ then $$\left(\frac{1}{n}\sum_{k=0}^{n-1}T^k(x_k)\right)_{n=1}^{\infty}\in L_p(\mathcal{M},\tau;c_0).$$ 

Furthermore, if $(x_k)_{k=0}^{\infty}\in L_p(\mathcal{M},\tau;\ell_\infty)$ and $x\in L_p(\mathcal{M},\tau)$ are such that $(x_k-x)_{k=0}^{\infty}\in L_p(\mathcal{M},\tau;c_0)$, then
$(\frac{1}{n}\sum_{k=0}^{n-1}T^k(x_k)-F_T(x))_{n=1}^{\infty}\in L_p(\mathcal{M},\tau;c_0),$ and so $\frac{1}{n}\sum_{k=0}^{n-1}T^k(x_k)\to F_T(x)$ b.a.u. as $n\to\infty$.
\end{teo}

\begin{proof}
First observe that if $(x_k)_{k=0}^{\infty}\in L_p(\mathcal{M},\tau;\ell_\infty)$ and $x\in L_p(\mathcal{M},\tau)$ are such that $(x_k-x)_{k=0}^{\infty}\in L_p(\mathcal{M},\tau;c_0)$, then
$$\left(\frac{1}{n}\sum_{k=0}^{n-1}T^k(x_k)-F_T(x)\right)_{n=1}^{\infty}=\left(\frac{1}{n}\sum_{k=0}^{n-1}T^k(x_k-x)\right)_{n=1}^{\infty}+\left(\frac{1}{n}\sum_{k=0}^{n-1}T^k(x)-F_T(x)\right)_{n=1}^{\infty}.$$ The second term on the right hand side is in $L_p(\mathcal{M},\tau;c_0)$ by \cite[Theorem 6.3]{jx}. Therefore the second claim will follow from the first since the first claim implies that $(\frac{1}{n}\sum_{k=0}^{n-1}T^k(x_k-x))_{n=1}^{\infty}\in L_p(\mathcal{M},\tau;c_0)$ since $(x_k-x)_{k=0}^{\infty}\in L_p(\mathcal{M},\tau;c_0)$. Hence, without loss of generality we will assume that $x=0$.

One can see from its definition and properties of the adjoint operator on $L_0(\mathcal{M},\tau)$ and the $\|\cdot\|_p$-norms that $(x_k)_{k=0}^{\infty}\in L_p(\mathcal{M},\tau;c_0)$ if and only if $(x_k^*)_{k=0}^{\infty}\in L_p(\mathcal{M},\tau;c_0)$. Since $L_p(\mathcal{M},\tau;c_0)$ is a subspace of $L_p(\mathcal{M},\tau;\ell_\infty)$, by proving the result for $\left(\frac{x_k+x_k^*}{2}\right)_{k=0}^{\infty}$ and $\left(\frac{x_k-x_k^*}{2i}\right)_{k=0}^{\infty}$ separately and then taking an appropriate linear combination without loss of generality we may assume that the sequence $(x_k)_{k=0}^{\infty}$ is a self-adjoint sequence in $L_p(\mathcal{M},\tau;c_0)$, i.e. we assume that $x_k^*=x_k$ for all $k\in\mathbb{N}_0$.

Fix $\epsilon>0$. Since $(x_k)_{k=0}^{\infty}\in L_p(\mathcal{M},\tau;c_0)$ is a self-adjoint sequence, it follows by \cite[Corollary 2.3]{lmz} that there exists $k^{\prime}\in\mathbb{N}$ and $c,d\in L_p(\mathcal{M},\tau)^+$ such that
$$\max\{\|c\|_p,\|d\|_p\}\leq\frac{\epsilon}{2C_p}\text{ and }-d\leq x_k\leq c\text{ for every }k\geq k^{\prime},$$ where $C_p$ denotes the strong type $(p,p)$ constant of the family $(M_n(T))_{n\in\mathbb{N}}$ found in \cite[Theorem 4.1]{jx}.

With $k^{\prime}$ now fixed there exists $m^{\prime}\in\mathbb{N}$ such that
$$\frac{1}{m^{\prime}}\left\|\sum_{k=0}^{k^{\prime}-1}T^k(x_k)\right\|_p\leq\frac{\epsilon}{2}.$$

If $n>k^{\prime}$, then since $T^k$ is positive for every $k\in\mathbb{N}_0$ it follows that
$$\frac{1}{n}\sum_{k=k^{\prime}}^{n-1}T^k(-d)\leq\frac{1}{n}\sum_{k=k^{\prime}}^{n-1}T^k(x_k)\leq\frac{1}{n}\sum_{k=k^{\prime}}^{n-1}T^k(c).$$ The left-most term here is larger than $-M_n(T)(d)$ and the right-most term is smaller that $M_n(T)(c)$ with respect to the ordering on $L_0(\mathcal{M},\tau)^h$ by the positivity of each $T^k$ and the fact that $c,d\geq0$. In other words,
$$-M_n(T)(d)\leq\frac{1}{n}\sum_{k=k^{\prime}}^{n-1}T^k(-d)\text{ and }\frac{1}{n}\sum_{k=k^{\prime}}^{n-1}T^k(c)\leq M_n(T)(c).$$ 

Since $c,d\in L_p(\mathcal{M},\tau)^{+}$ and since $(M_n(T))_{n\in\mathbb{N}}$ is of strong type $(p,p)$ with constant $C_p$, it follows that there exists $f,g\in L_p(\mathcal{M},\tau)^+$ such that $\|f\|_{p}\leq C_p\|c\|_p$, $\|g\|_p\leq C_p\|d\|_p$, and
$$0\leq M_n(T)(c)\leq f\text{ and }0\leq M_n(T)(d)\leq g\text{ for every }n\in\mathbb{N}.$$
It follows from all of this that
$$-g\leq -M_n(T)(d)\leq
\frac{1}{n}\sum_{k=k^{\prime}}^{n-1}T^k(x_k)
\leq M_n(T)(c)\leq f
\text{ for every }n\geq k^{\prime}.$$

Now for every $n\geq\max\{m^{\prime},k^{\prime}\}$ it holds that
$$%\begin{align*}
\frac{1}{n}\sum_{k=0}^{n-1}T^k(x_k)
=\frac{1}{n}\sum_{k=0}^{k^{\prime}-1}T^k(x_k)+\frac{1}{n}\sum_{k=k^{\prime}}^{n-1}T^k(x_k) \\
%\leq\frac{1}{n}\left|\sum_{k=0}^{k_0-1}T^k(x_k)\right|+f 
\leq\frac{1}{m^{\prime}}\left|\sum_{k=0}^{k^{\prime}-1}T^k(x_k)\right|+f
$$%\end{align*}
and 
\begin{align*}
	\left\|\frac{1}{m^{\prime}}\left|\sum_{k=0}^{k^{\prime}-1}T^k(x_k)\right|+f\right\|_p
	\leq&\left\|\frac{1}{m^{\prime}}\left|\sum_{k=0}^{k^{\prime}-1}T^k(x_k)\right|\right\|_p+\|f\|_p \\
	\leq&\frac{1}{m^{\prime}}\left\|\sum_{k=0}^{k^{\prime}-1}T^k(x_k)\right\|_p+C_p\|c\|_p \\
	\leq&\frac{\epsilon}{2}+C_p\cdot\frac{\epsilon}{2C_p}
	=\epsilon.
\end{align*}
A similar argument for the same range of $n$ shows that
$$\frac{1}{n}\sum_{k=0}^{n-1}T^k(x_k)\geq-\left(\frac{1}{m^{\prime}}\left|\sum_{k=0}^{k^{\prime}-1}T^k(x_k)\right|+g\right)$$ and $$\left\|\frac{1}{m^{\prime}}\left|\sum_{k=0}^{k^{\prime}-1}T^k(x_k)\right|+g\right\|_p\leq\epsilon.$$ 

Since $\epsilon>0$ was arbitrary, it follows that $\left(\frac{1}{n}\sum_{k=0}^{n-1}T^k(x_k)\right)_{n=1}^{\infty}\in L_p(\mathcal{M},\tau;c_0)$ by \cite[Corollary 2.3]{lmz}, finishing the proof.
\end{proof}

\begin{rem}\label{GeneralvonNeumannAlgebras}
The main strategy for the proof of Theorem \ref{GeneralBreimanErgodicTheorem} mostly relied on properties of the spaces $L_p(\mathcal{M},\tau;\ell_\infty)$ and $L_p(\mathcal{M},\tau;c_0)$. Furthermore, all of the properties used there all generalize to the general (not necessarily semifinite) von Neumann algebra setup as mentioned in \cite[Sections 2 and 7]{jx} and \cite[Remark 5.6]{lmz}. As such, Theorem \ref{GeneralBreimanErgodicTheorem} extends to the noncommutative Haagerup spaces $L_p(\mathcal{N})=L_p(\mathcal{N},\sigma)$ associated to a von Neumann algebra $\mathcal{N}$ with a normal semifinite faithful weight $\sigma:\mathcal{N}^{+}\to[0,\infty]$ when $T:\mathcal{N}\to\mathcal{N}$ is a completely positive contraction that commutes with the modular automorphism group of $\sigma$ and is such that $\sigma\circ T\leq \sigma$, writing $T$ for its continuous extensions to $L_p(\mathcal{N})$ as well. We leave the details to the reader.
\end{rem}

Using the alternative characterization of self-adjoint elements of $L_p(\mathcal{M},\tau;\ell_\infty)$ one can easily modify the previous proof to show that $(x_k)_{k=0}^{\infty}\in L_p(\mathcal{M},\tau;\ell_\infty)$ implies that $(\frac{1}{n}\sum_{k=0}^{n-1}T^k(x_k))_{n=1}^{\infty}\in L_p(\mathcal{M},\tau;\ell_\infty)$. Furthermore, that proof itself primarily shows that $L_p(\mathcal{M},\tau;c_0)$ is invariant under this assignment. Due to these, one can use the closed graph theorem to show the following. We leave the details to the reader.

\begin{cor}\label{Corollary}
Let $T\in DS^+(\mathcal{M},\tau)$ and $1<p<\infty$. For each $(x_k)_{k=0}^{\infty}\in L_p(\mathcal{M},\tau;\ell_\infty)$ define
$$M\left((x_k)_{k=0}^{\infty}\right):=\left(\frac{1}{n}\sum_{k=0}^{n-1}T^k(x_k)\right)_{n=1}^{\infty}.$$ Then $M$ defines a bounded linear operator from $L_p(\mathcal{M},\tau;\ell_\infty)$ to itself with a norm of at most $6C_p$ which leaves $L_p(\mathcal{M},\tau;c_0)$ invariant, where $C_p$ denotes the strong type $(p,p)$ constant of the family $(M_n(T))_{n\in\mathbb{N}}$.
\end{cor}

We will now prove the direct extension of Breiman's individual ergodic theorem to the von Neumann algebra setting. As previously mentioned, the proof will be similar to that of Theorem \ref{GeneralBreimanErgodicTheorem}, but the lack of a strong type $(1,1)$ maximal inequality for $(M_n(T))_{n\in\mathbb{N}}$ means a bit more care will be required in the proof in showing both that b.a.u. convergence will occur and what the limit must be. 

By requiring $(\mathcal{M},\tau)$ to be tracial we will be able to use the fact that $M_n(T)(x)\to\tau(x)\textbf{1}$ b.a.u. as $n\to\infty$ for every $x\in L_1(\mathcal{M},\tau)$ when $T$ is induced by an ergodic normal $\tau$-preserving $*$-homomorphism of $\mathcal{M}$. This particular form of the limit plays an important role in the proof. In particular, if $F_T(x)$ denotes the b.a.u. limit for $x\in L_1(\mathcal{M},\tau)$ (which exists by \cite[Theorem 2]{ye}), then part of the proof requires finding an upper bound on $\|eF_T(c)e\|_\infty$ for some $c\in L_1(\mathcal{M},\tau)^{+}$ and $e\in\mathcal{P(M)}$, which has no easy bounds even if $\|c\|_1$ is small. However, if $F_T(c)=\tau(c)\textbf{1}$, then since $\tau(c)=\|c\|_1$ we would find that $\|eF_T(c)e\|_\infty=\tau(c)\|e\|_\infty\leq\|c\|_1$, which is then used for an upper bound on $\|eM_n(T)(c)e\|_\infty$ in terms of $\|c\|_1$.

The identification of the actual limit operator follows from the facts that the map $F_T:L_1(\mathcal{M},\tau)\to L_1(\mathcal{M},\tau)$ defined by $F_T(x):=\|\cdot\|_1-\lim_{n\to\infty}M_n(T)(x)$ for $x\in L_1(\mathcal{M},\tau)$ is continuous when $p=1$ since $\tau(\textbf{1})<\infty$ by \cite[Lemma 2.4]{comli1} and the mean ergodic theorem, that $F_T(x)=\tau(x)\textbf{1}$ for $x\in L_2(\mathcal{M},\tau)$ by \cite[Proposition 5.1]{li2} since b.a.u. and $L_1(\mathcal{M},\tau)$-norm and $L_2(\mathcal{M},\tau)$-norm limits agree because the measure topology on $L_0(\mathcal{M},\tau)$ is Hausdorff, and since $L_2(\mathcal{M},\tau)$ is dense in $L_1(\mathcal{M},\tau)$.

%\vfill

\begin{teo}\label{DirectBreimanErgodicTheorem}
Let $(\mathcal{M},\tau)$ be a tracial von Neumann algebra and $T:L_1(\mathcal{M},\tau)\to L_1(\mathcal{M},\tau)$ be induced by an ergodic normal $\tau$-preserving $*$-homomorphism of $\mathcal{M}$. If $(x_k)_{k=0}^{\infty}\in L_1(\mathcal{M},\tau;\ell_\infty)$ and $x\in L_1(\mathcal{M},\tau)$ are such that $(x_k-x)_{k=0}^{\infty}\in L_1(\mathcal{M},\tau;c_0)$, then
$$\frac{1}{n}\sum_{k=0}^{n-1}T^k(x_k)\to\tau(x)\textbf{1}\text{ b.a.u. as }n\to\infty.$$
\end{teo}

\begin{proof}
As in the proof of Theorem \ref{GeneralBreimanErgodicTheorem} it suffices to prove the case where $x=0$ and $(x_k)_{k=0}^{\infty}\in L_1(\mathcal{M},\tau;c_0)$ is such that $x_k^*=x_k$ for every $k\in\mathbb{N}_0$.

Let $\epsilon,\delta>0$. By \cite[Corollary 2.3]{lmz} there exists $k^{\prime}\in\mathbb{N}_0$ and $c,d\in L_1(\mathcal{M},\tau)^+$ such that
$$\max\{\|c\|_1,\|d\|_1\}\leq\frac{\delta}{24}\text{ and }-d\leq x_k\leq c\text{ for every }k\geq k^{\prime}.$$
For the same reason as in the proof of Theorem \ref{GeneralBreimanErgodicTheorem} it follows that
$$-M_n(T)(d)\leq\frac{1}{n}\sum_{k=k^{\prime}}^{n-1}T^k(x_k)\leq M_n(T)(c)$$ in $L_1(\mathcal{M},\tau)^{h}$ for every $n>k^{\prime}$.

Let $e\in\mathcal{P(M)}$ be such that $\tau(e^\perp)\leq\epsilon$, $T^k(x_k)e,T^k(c)e,T^k(d)e\in\mathcal{M}$ for every $k\in\mathbb{N}_0$, and both $\|e(M_n(T)(c)-\tau(c)\textbf{1})e\|_\infty\to0$ and $\|e(M_n(T)(d)-\tau(d)\textbf{1})e\|_\infty\to0$ as $n\to\infty$. Such a projection $e$ exists by taking the intersection of the projections one obtains from the definition of each $T^k(x_k)$, $T^k(c)$, and $T^k(d)$ being $\tau$-measurable and the fact that $(M_n(T)(z))_{n=1}^{\infty}$ converges to $\tau(z)\textbf{1}$ b.a.u. for every $z\in L_1(\mathcal{M},\tau)$ due to the assumptions on $\mathcal{M}$ and $T$ by \cite[Theorem 2]{ye}.

Let $A_n:=\frac{1}{n}\sum_{k=k^{\prime}}^{n-1}T^k(x_k)$ for every $n> k^{\prime}$. Then the above implies that
$$-M_n(T)(d)\leq A_n\leq M_n(T)(c)$$ in $L_1(\mathcal{M},\tau)^{h}$, which then further implies that
$$-eM_n(T)(d)e\leq eA_ne\leq eM_n(T)(c)e$$ in $\mathcal{M}^{h}$ for every such $n$. Hence, it follows that
\begin{align}\label{LimSupArgument}
0
\leq\|e(A_n+M_n(T)(d))e\|_\infty,\|e(M_n(T)(c)-A_n)e\|_\infty
\leq\|eM_n(T)(c+d)e\|_\infty
\end{align} 
for every $n>k^{\prime}$.

Since $k^{\prime}$ is fixed there exists $m^{\prime}\in\mathbb{N}$ such that
$$\frac{1}{m^{\prime}}\left\|\sum_{k=0}^{k^{\prime}-1}eT^k(x_k)e\right\|_\infty\leq\frac{\delta}{4}.$$ 

Let $n^{\prime}\in\mathbb{N}$ be such that $n\geq n^{\prime}$ implies that
\begin{align*}
	\|e(M_n(T)(c)-\tau(c)\textbf{1})e\|_\infty
	\leq&|\tau(c+d)| \\
	\|e(M_n(T)(d)-\tau(d)\textbf{1})e\|_\infty
	\leq&|\tau(c+d)| \\
	\|eM_n(T)(c+d)e\|_\infty
	\leq&\|e(\tau(c+d)\textbf{1})e\|_\infty+\frac{\delta}{24}
\end{align*}
Such an $n^{\prime}$ exists from the fact that $\|e(M_n(T)(z)-\tau(z)\textbf{1})e\|_\infty\to0$ as $n\to\infty$ for each $z\in\{c,d,c+d\}$ and the reverse triangle inequality.

Since $c,d\in L_1(\mathcal{M},\tau)^{+}$ one can find that $\tau(c+d)=\tau(|c+d|)=\|c+d\|_1$, so that $\|e(\tau(c+d)\textbf{1})e\|_\infty=|\tau(c+d)|\|e\|_\infty\leq\|c+d\|_1$ since $\|e\|_\infty\leq1$. From all this it follows that if $m,n\geq\max\{k^{\prime},m^{\prime},n^{\prime}\}$, then all of this implies that
\begin{align*}
	\|e(A_n-A_m)e\|_\infty
	\leq&\|e(A_n+M_n(T)(d))e\|_\infty+\|e(M_n(T)(d)+M_m(T)(c))e\|_\infty \\
	& \ \ \ \ \ +\|e(M_m(T)(c)-A_m)e\|_\infty \\
	\leq&\|eM_n(T)(c+d)e\|_\infty+\|e(M_n(T)(d)+M_m(T)(c))e\|_\infty \\
	& \ \ \ \ \ +\|eM_m(T)(c+d)e\|_\infty \\
	\leq&\|eM_n(T)(c+d)e\|_\infty+\|e(M_n(T)(d)-\tau(d)\textbf{1})e\|_\infty \\
	&  \ \ \ \ \ +\|e(M_m(T)(c)-\tau(c)\textbf{1})e\|_\infty +\|e(\tau(c+d)\textbf{1})e\|_\infty \\
	& \ \ \ \ \ +\|eM_m(T)(c+d)e\|_\infty \\
	\leq&\left(\|\tau(c+d)e\|_\infty+\frac{\delta}{24}\right)+|\tau(c+d)|+|\tau(c+d)| \\
	& \ \ \ \ \ +\|\tau(c+d)e\|_\infty+\left(\|\tau(c+d)e\|_\infty+\frac{\delta}{24}\right) \\
	%\leq&5\tau(|c+d|)+\frac{2\delta}{24} \\
	%\leq&5\tau(|c+d|)+\frac{\delta}{12}
	\leq&5\|c+d\|_1+\frac{\delta}{12} \\
	\leq&5(\|c\|_1+\|d\|_1)+\frac{\delta}{12}
	\leq5\cdot\frac{2\delta}{24}+\frac{2\delta}{24}
	=\frac{\delta}{2},
\end{align*}
from which it follows that
\begin{align*}
	\bigg\|e\bigg(\frac{1}{n}\sum_{k=0}^{n-1}&T^k(x_k)-\frac{1}{m}\sum_{k=0}^{m-1}T^k(x_k)\bigg)e\bigg\|_\infty \\
	\leq&\frac{1}{n}\left\|\sum_{k=0}^{k^{\prime}-1}eT^k(x_k)e\right\|_\infty +\frac{1}{m}\left\|\sum_{k=0}^{k^{\prime}-1}eT^k(x_k)e\right\|_\infty+\|e(A_n-A_m)e\|_\infty \\
	\leq&\frac{\delta}{4}+\frac{\delta}{4}+\frac{\delta}{2}=\delta.
\end{align*}
Since $m,n\geq\max\{k^{\prime},m^{\prime},n^{\prime}\}$, $\delta>0$, and $\epsilon>0$ were all arbitrary, it follows that the sequence $(\frac{1}{n}\sum_{k=0}^{n-1}T^k(x_k))_{n=1}^{\infty}$ is b.a.u. Cauchy, and so it converges to some $\widehat{x}\in L_0(\mathcal{M},\tau)$ b.a.u. as $n\to\infty$ by \cite[Theorem 2.3]{cls}. 

We will now show that $\widehat{x}=0$. To do this, we will find a subsequence of $(\frac{1}{n}\sum_{k=0}^{n-1}T^k(x_k))_{k=1}^{\infty}$ inductively which converges to $0$ with respect to the norm of $L_1(\mathcal{M},\tau)$. To start, define $n_0:=1$.

Fix $j\in\mathbb{N}$, and let $k_j\in\mathbb{N}_0$ and $c_j,d_j\in L_1(\mathcal{M},\tau)^+$ be such that
$$\max\{\|c_j\|_1\|d_j\|_1\}\leq\frac{1}{6j}\text{ and }-d_j\leq x_k\leq c_j\text{ for every }k\geq k_j.$$
Note that this implies that $0\leq x_k+d_j\leq c_j+d_j$, and so $\|x_k+d_j\|_1\leq\|c_j+d_j\|_1$, for every $k\geq k_j$.

Let $m_j\in\mathbb{N}$ be such that
$$\frac{1}{m_j}\left\|\sum_{k=0}^{k_j-1}T^k(x_k)\right\|_1\leq\frac{1}{2j}.$$
If $n_{j}:=\max\{n_{j-1},k_j,m_j\}+1$, then we find that
\begin{align*}
	\left\|\frac{1}{n_j}\sum_{k=0}^{n_{j}-1}T^k(x_k)\right\|_1
	\leq&\frac{1}{n_j}\left\|\sum_{k=0}^{k_{j}-1}T^k(x_k)\right\|_1+\frac{1}{n_j}\sum_{k=k_{j}}^{n_{j}-1}\|T^k(x_k)\|_1 \\
	\leq&\frac{1}{m_j}\left\|\sum_{k=0}^{k_j-1}T^k(x_k)\right\|_1+\frac{1}{n_j}\sum_{k=k_j}^{n_j-1}\|x_k\|_1 \\
	\leq&\frac{1}{2j}+\frac{1}{n_j}\sum_{k=k_j}^{n_j-1}\bigg[\|x_k+d\|_1+\|d\|_1\bigg] \\
	\leq&\frac{1}{2j}+\frac{1}{n_j}\sum_{k=0}^{n_j-1}\bigg[\|c+d\|_1+\|d\|_1\bigg] \\
	\leq&\frac{1}{2j}+\bigg[\|c\|_1+2\|d\|_1\bigg]
	\leq\frac{1}{2j}+\frac{1}{6j}+\frac{2}{6j}=\frac{1}{j}.
\end{align*}
Since $j\in\mathbb{N}$ was arbitrary, and since $n_j\leq n_{j+1}$ for every $j\in\mathbb{N}$, it follows that $(\frac{1}{n_j}\sum_{k=0}^{n_j-1}T^k(x_k))_{j=0}^{\infty}$ is a subsequence of $(\frac{1}{n}\sum_{k=0}^{n-1}T^k(x_k))_{n=1}^{\infty}$ which converges to $0$ with respect to the norm of $L_1(\mathcal{M},\tau)$. Since the same subsequence would have the same b.a.u. limit as the original sequence since both b.a.u. and $L_1(\mathcal{M},\tau)$-norm limits agree because both convergences imply convergence in measure and $L_0(\mathcal{M},\tau)$ is Hausdorff with respect to the measure topology, it follows that the b.a.u. limit of $(\frac{1}{n}\sum_{k=0}^{n-1}T^k(x_k))_{n=1}^{\infty}$ must be equal to $0$.
\end{proof}

We now take a short interlude to compare the theory of the noncommutative vector-valued $L_p$-spaces with the commutative ones. We do this to justify the assumptions made on $(x_k)_{k=0}^{\infty}$ made in the results above since what we use does not immediately resemble what was stated in the original result of Breiman \cite{Breiman57} (which we stated as Theorem \ref{BreimanOriginalResult} in Section \ref{Intro}). In particular, we will show that this formulation is implied by what was made by Breiman. Interestingly, a large part of Breiman's proof %of \cite[Theorem 1]{Breiman57} 
is dedicated to showing that this condition holds and what it implies. We use \cite{defant} as our primary reference for this comparison.

Let $(\Omega,\mathcal{F},\nu)$ be a $\sigma$-finite measure space and fix $1\leq p<\infty$. Define $L_p(\Omega,\mathcal{F},\nu;\ell_\infty)$ to be the set of all sequences $f=(f_k)_{k=0}^{\infty}\subset L_p(\Omega,\mathcal{F},\nu)$ whose maximal function $f^*$ is an element of $L_p(\Omega,\mathcal{F},\nu)$, where
$$f^*(\omega):=\sup_{k\in\mathbb{N}_0}|f_k(\omega)|\text{ for every }\omega\in\Omega.$$ 
Equipped with the norm
$$\|(f_k)_{k=0}^{\infty}\|_{L_p(\Omega,\mathcal{F},\nu;\ell_\infty)}:=\|f^*\|_p$$ it is known that $L_p(\Omega,\mathcal{F},\nu;\ell_\infty)$ is a Banach space with coordinate-wise operations. If $|f|:=(|f_k|)_{k=0}^\infty$ (where $|g|(\omega):=|g(\omega)|$ for $\omega\in\Omega$ and $g\in L_p(\Omega,\mathcal{F},\nu)$) then it follows that $f^*=|f|^*$, so that $(f_k)_{k=0}^{\infty}\in L_p(\Omega,\mathcal{F},\nu;\ell_\infty)$ if and only if $(f_k)_{k=0}^{\infty}\in L_p(\Omega,\mathcal{F},\nu;\ell_\infty)$, which is a feature that does not hold in the noncommutative setting.

It was shown in \cite[Lemma 1]{defant} that if $(f_k)_{k=0}^{\infty}\in L_p(\Omega,\mathcal{F},\nu)$, then $(f_k)_{k=0}^{\infty}\in L_p(\Omega,\mathcal{F},\nu;\ell_\infty)$ if and only if there exists $(z_k)_{k=0}^{\infty}\subset L_\infty(\Omega,\mathcal{F},\nu)$ with $\sup_{k\in\mathbb{N}_0}\|z_k\|_\infty<\infty$ and $g\in L_p(\Omega,\mathcal{F},\nu)$ such that
$$f_k=z_kg\text{ for every }k\in\mathbb{N}_0.$$
Furthermore,
$$\|(f_k)_{k=0}^{\infty}\|_{L_p(\Omega,\mathcal{F},\nu;\ell_\infty)}=\inf_{f_k=z_kg}\sup_{k\in\mathbb{N}_0}\|z_k\|_\infty\|g\|_p,$$ where the infimum runs over all factorizations of $(f_k)_{k=0}^{\infty}$ mentioned above. 

For every $z\in L_\infty(\Omega,\mathcal{F},\nu)$ define the operator $M_z:L_2(\Omega,\mathcal{F},\nu)\to L_2(\Omega,\mathcal{F},\nu)$ defined by
$$M_zf(\omega):=z(\omega)f(\omega)\text{ for every }\omega\in\Omega\text{ and }f\in L_2(\Omega,\mathcal{F},\nu).$$ Let $\mathcal{Z}$ denote the set of all such operators, and let $\tau_{\nu}$ be the trace on $\mathcal{Z}$ defined by $\tau_{\nu}(M_z):=\int_{\Omega}z\ d\nu$. It is known that $\mathcal{Z}$ is a von Neumann algebra on $L_2(\Omega,\mathcal{F},\nu)$.

Let $g\in L_p(\Omega,\mathcal{F},\nu)$ and define $\mathcal{D}_{g}:=\{f\in L_2(\Omega,\mathcal{F},\nu):gf\in L_2(\Omega,\mathcal{F},\nu)\}$. Define $M_g:\mathcal{D}_{g}\to L_2(\Omega,\mathcal{F},\nu)$ by 
$$M_{g}f:=g(\omega)f(\omega)\text{ for every }\omega\in\Omega\text{ and every }f\in\mathcal{D}_{g}.$$ Then $M_g\in L_p(\mathcal{Z},\tau_{\nu})$ with $\|M_g\|_p=\|g\|_p$. Furthermore, to every $T\in L_p(\mathcal{Z},\tau_{\nu})$ there exists $g\in L_p(\Omega,\mathcal{F},\nu)$ such that $T=M_g$.

From their definitions, the above tells us that if $(f_k)_{k=0}^{\infty}\in L_p(\Omega,\mathcal{F},\nu)$, then
$$(f_k)_{k=0}^{\infty}\in L_p(\Omega,\mathcal{F},\nu;\ell_\infty)\text{ if and only if } (M_{f_k})_{k=0}^{\infty}\in L_p(\mathcal{Z},\tau_{\nu};\ell_\infty)$$
and
$$\left\|{\sup_{k\in\mathbb{N}_0}}^{+}M_{f_k}\right\|_p=\left\|\sup_{k\in\mathbb{N}_0}|f_k|\right\|_p,$$ where we emphasize that the expression on the left is the norm of the sequence $(M_{f_k})_{k=0}^{\infty}\in L_p(\mathcal{Z},\tau_{\nu};\ell_\infty)$ and the quantity on the right is the norm of $(f_k)_{k=0}^{\infty}\in L_p(\Omega,\mathcal{F},\nu;\ell_\infty)$.

The space $L_p(\Omega,\mathcal{F},\nu;c_0)$ is defined similarly, and coincides similarly with the space $L_p(\mathcal{M},\tau_{\nu};c_0)$.

It is common to identify $f$ with $M_f$ for every $f\in L_p(\Omega,\mathcal{F},\nu)$, so the left hand side of the above would be written as $\|{\sup_{k\in\mathbb{N}_0}}^{+}f_k\|_p$. We will only need this identification for the proof of the next result, since we wish to discuss sequences in $L_p(\Omega,\mathcal{F},\nu;\ell_\infty)$ and $L_p(\Omega,\mathcal{F},\nu;c_0)$ while using properties of the spaces $L_p(\mathcal{Z},\tau_{\nu};\ell_\infty)$ and $L_p(\mathcal{Z},\tau_{\nu};c_0)$ which were proven in \cite{lmz}.

%where the expression on the left is the norm of $(f_k)_{k=0}^{\infty}$ as an element of the space $L_p(L_\infty(\Omega,\mathcal{F},\mu),\int_{\Omega}\cdot d\mu;\ell_\infty)$ (associated to the von Neumann algebra $L_\infty(\Omega,\mathcal{F},\mu)$ acting on $L_2(\Omega,\mathcal{F},\mu)$) and the expression on the right is the $L_p$-norm of the maximal function $\sup_{k\in\mathbb{N}_0}|f_k|$. Due to this, if $\|\sup_{k\in\mathbb{N}_0}|f_k|\|_p<\infty$ then it can be shown that $(f_k)_{k=0}^{\infty}\in L_p(L_\infty(\Omega,\mathcal{F},\mu),\int_{\Omega}\cdot d\mu;\ell_\infty)$.

\begin{lm}\label{WhyAssumption}
	Let $(\Omega,\mathcal{F},\nu)$ be a $\sigma$-finite measure space. Assume that $(g_k)_{k=0}^{\infty}\subset L_1(\Omega,\mathcal{F},\nu)$ is such that $\int_{\Omega}\sup_{k\in\mathbb{N}_0}|g_k(\omega)|d\nu(\omega)<\infty$, and assume that there exists a measurable function $g:\Omega\to\mathbb{C}$ such that $g_k\to g$ $\nu$-a.e. Then it holds that $(g_k)_{k=0}^{\infty}\in L_1(\Omega,\mathcal{F},\nu;\ell_\infty)$, $g\in L_1(\Omega,\mathcal{F},\nu)$, and $$(g_k-g)_{k=0}^{\infty}\in L_1(\Omega,\mathcal{F},\nu;c_0).$$
\end{lm}
\begin{proof}
	Define $g^*(\omega):=\sup_{k\in\mathbb{N}_0}|g_k(\omega)|$, noting that $g^*\in L_1(\Omega,\mathcal{F},\nu)$ by assumption; hence $(g_k)_{k=0}^{\infty}\in L_1(\Omega,\mathcal{F},\nu;\ell_\infty)$. Since $(g)_{k=0}^{\infty}\in L_1(\Omega,\mathcal{F},\nu;\ell_\infty)$ as well it follows that $(g_k-g)_{k=0}^{\infty}=(g_k)_{k=0}^{\infty}-(g)_{k=0}^{\infty}$, which then implies that $(|g_k-g|)_{k=0}^{\infty}\in L_p(\Omega,\mathcal{F},\nu;\ell_\infty)$.
	
	For every $N\in\mathbb{N}_0$ define $G_N(\omega):=\sup_{k\geq N}|g_k(\omega)-g(\omega)|$ for $\omega\in\Omega$. From their definition it holds that $G_N\geq G_{N+1}$ as functions on $\Omega$ for every $N\in\mathbb{N}_0$.
	
	By the dominated convergence theorem it follows from $g_k\to g$ $\nu$-a.e., $|g_k|\leq g^*$ for every $k\in\mathbb{N}_0$, and $g^*\in L_1(\Omega,\mathcal{F},\nu)$ that $g\in L_1(\Omega,\mathcal{F},\nu)$.
	
	Since $g_k\to g$ $\nu$-a.e., there exists a set $\Omega^{\prime}\in\mathcal{F}$ with $\nu(\Omega\setminus\Omega^{\prime})=0$ such that $\omega\in\Omega^{\prime}$ implies $g_k(\omega)\to g(\omega)$ as $k\to\infty$.
	
	Fix $\omega\in\Omega^{\prime}$ and $\epsilon>0$. Then there exists $K\in\mathbb{N}_0$ such that $k\geq K$ implies that $|g_k(\omega)-g(\omega)|<\epsilon$. From this we find that $N\geq K$ implies that
	$$0\leq G_N(\omega)=\sup_{k\geq N}|g_k(\omega)-g(\omega)|\leq\sup_{k\geq K}|g_k(\omega)-g(\omega)|\leq\epsilon.$$ Since $\epsilon>0$ was arbitrary it follows that $G_N(\omega)\to0$ as $N\to\infty$, and since $\omega\in\Omega^{\prime}$ was arbitrary it follows that $G_N\to0$ $\mu$-a.e.
	
	Since $0\leq G_N\leq G_0\leq g^*+|g|$ for every $N\in\mathbb{N}_0$, and since $g,g^*\in L_1(\Omega,\mathcal{F},\nu)$ implies that $g^*+|g|\in L_1$, it follows that $G_N\in L_1(\Omega,\mathcal{F},\nu)$ for every $N\geq0$, and so the dominated convergence theorem implies
	$$\lim_{N\to\infty}\left\|\sup_{k\geq N}|g_k-g|\right\|_1
	=\lim_{N\to\infty}\|G_N\|_1
	=\lim_{N\to\infty}\int_{\Omega}G_N\,d\mu
	=\int_{\Omega}\lim_{N\to\infty}G_N\,d\mu
	%=\int_{\Omega}0d\mu
	=0$$ since $G_N\to0$ $\nu$-a.e. as $N\to\infty$. 
	Since $(|g_k-g|)_{k=0}^{\infty}$ is already an element of $L_p(\Omega,\mathcal{F},\nu;\ell_\infty)$ it follows by \cite[Lemma 2.2]{lmz} that $(|g_k-g|)_{k=0}^{\infty}\in L_1(\Omega,\mathcal{F},\nu;c_0)$. 
	By definition there then exists $f\in L_1(\Omega,\mathcal{F},\nu)$ and $(z_k)_{k=0}^{\infty}\in L_\infty(\Omega,\mathcal{F},\nu)$ with $\lim_{k\to\infty}\|z_k\|_\infty=0$ such that $$|g_k-g|=z_kf\text{ for every }k\in\mathbb{N}_0.$$ By the polar decomposition of functions, for every $k\in\mathbb{N}_0$ there exists $u_k\in L_\infty(\Omega,\mathcal{F},\nu)$ with $\|u_k\|_\infty\leq1$ such that $g_k-g=u_k|g_k-g|$. Hence
	$$g_k-g=u_k|g_k-g|=(u_kz_k)f\text{ for every }k\in\mathbb{N}_0$$ with $$\lim_{k\to\infty}\|u_kz_k\|_\infty\leq\lim_{k\to\infty}\|u_k\|_\infty\|z_k\|_\infty\leq\lim_{k\to\infty}\|z_k\|_\infty=0,$$ which implies $\lim_{k\to\infty}\|u_kz_k\|_\infty=0$, and so $(g_k-g)_{k=0}^{\infty}\in L_1(\Omega,\mathcal{F},\nu;c_0)$.
\end{proof}

\section{Further Generalizations of Breiman's Ergodic Theorem}\label{SupplementaryResults}

In this section we will prove a number of variations of the ergodic theorems considered in the previous section. In particular, we will consider weighted and subsequential versions of the result, most of which follow nearly immediately from their unmodified counterparts. In Theorem \ref{EntangledAverages2} we consider averages of the form
$$\frac{1}{n}\sum_{k=0}^{n-1}T^kRS^k$$ resembling the entangled ergodic averages considered in \cite{ekk,lieb}; in our version of these results, we assume that $R:L_p(\mathcal{M},\tau)\to L_p(\mathcal{M},\tau)$ is a positive bounded linear operator and $S,T\in DS^+(\mathcal{M},\tau)$ are such that $S^d$ is positive as a Hilbert space operator on $L_2(\mathcal{M},\tau)$ for some $d\in\mathbb{N}$. Finally, we relate the results of this article to the results about $T$-admissible processes considered in \cite{ob5}.

\subsection{Weighted Averages}

In this subsection we will study weighted versions of the averages studied above. We start with a construction of finding new sequences in $L_p(\mathcal{M},\tau;c_0)$ from old ones.

\begin{pro}\label{BoundedWeightExample}
If $\alpha=(\alpha_k)_{k=0}^{\infty}\in \ell_\infty$ is fixed, then the mapping 
\begin{align*}
L_p(\mathcal{M},\tau;\ell_\infty)&\to L_p(\mathcal{M},\tau;\ell_\infty)\\
(x_k)_{k=0}^{\infty}&\mapsto(\alpha_k x_k)_{k=0}^{\infty}
\end{align*}
is linear and bounded with
$$\left\|{\sup_{k\in\mathbb{N}_0}}^{+}\alpha_k x_k\right\|_p\leq\|\alpha\|_{\ell_\infty}\left\|{\sup_{k\in\mathbb{N}_0}}^{+}x_k\right\|_p.$$	

\noindent Furthermore, in either of the following cases
\begin{itemize}
\item[(i)] $(x_k)_{k=0}^{\infty}\in L_p(\mathcal{M},\tau;c_0)$ and $(\alpha_k)_{k=0}^{\infty}\in\ell_\infty$, 
\item[(ii)] $(x_k)_{k=0}^{\infty}\in L_p(\mathcal{M},\tau;\ell_\infty)$ and $(\alpha_k)_{k=0}^{\infty}\in c_0$,
\end{itemize}
it holds that $(\alpha_k x_k)_{k=0}^{\infty}\in L_p(\mathcal{M},\tau;c_0)$.
\end{pro}
\begin{proof}		
Fix $\epsilon>0$, and consider a factorization of $(x_k)_{k=0}^{\infty}\in L_p(\mathcal{M},\tau;\ell_\infty)$ as
$$x_k=ay_kb\text{ for every }k\in\mathbb{N}_0,$$ where $a,b\in L_{2p}(\mathcal{M},\tau)$ and $(y_k)_{k=0}^{\infty}\subset\mathcal{M}$ are such that $\sup_{k\in\mathbb{N}_0}\|y_k\|_\infty<\infty$ and $\|a\|_{2p}\|b\|_{2p}\sup_{k\in\mathbb{N}_0}\|y_k\|_\infty\leq\|{\sup_{k\in\mathbb{N}_0}}^{+}x_k\|_p+\epsilon$. Then we can factor $(\alpha_k x_k)_{k=0}^{\infty}$ as
$$\alpha_k x_k=a(\alpha_k y_k)b\text{ for every }k\in\mathbb{N}_0,$$ noting that $\sup_{k\in\mathbb{N}_0}\|\alpha_k y_k\|_\infty=\sup_{k\in\mathbb{N}_0}|\alpha_k|\|y_k\|_\infty\leq\|\alpha\|_{\ell_\infty}\sup_{k\in\mathbb{N}_0}\|y_k\|_\infty<\infty$. Then, since  
\begin{align*}
\|a\|_{2p}\|b\|_{2p}\sup_{k\in\mathbb{N}_0}\|\alpha_k y_k\|_\infty
\leq&\|\alpha\|_{\ell_\infty}\bigg(\|a\|_{2p}\|b\|_{2p}\sup_{k\in\mathbb{N}_0}\|y_k\|_\infty\bigg) \\
\leq&\|\alpha\|_{\ell_\infty}\left(\left\|{\sup_{k\in\mathbb{N}_0}}^{+}x_k\right\|_p+\epsilon\right)
,
\end{align*} 
by taking the infimum over all appropriate factorizations of $(\alpha_k x_k)_{k=0}^{\infty}$ above it follows that
$$\left\|{\sup_{k\in\mathbb{N}_0}}^{+}\alpha_k x_k\right\|_p\leq\|\alpha\|_{\ell_\infty}\left(\left\|{\sup_{k\in\mathbb{N}_0}}^{+}x_k\right\|_p+\epsilon\right).$$ Since $\epsilon>0$ was arbitrary this implies that
$$\left\|{\sup_{k\in\mathbb{N}_0}}^{+}\alpha_k x_k\right\|_p\leq\|\alpha\|_{\ell_\infty}\left\|{\sup_{k\in\mathbb{N}_0}}^{+}x_k\right\|_p,$$ finishing the proof of the first part.

Now assume that $(x_k)_{k=0}^{\infty}\in L_p(\mathcal{M},\tau;c_0)$ and consider a factorization $x_k=ay_kb$ with $a,b\in L_{2p}(\mathcal{M},\tau)$ and $(y_k)_{k=0}^{\infty}\subset\mathcal{M}$ such that  $\lim_{k\to\infty}\|y_k\|_\infty=0$. Then
$$\lim_{k\to\infty}\|\alpha_k y_k\|_\infty=\lim_{k\to\infty}|\alpha_k|\|y_k\|_\infty\leq\|\alpha\|_{\ell_\infty}\lim_{k\to\infty}\|y_k\|_\infty=0,$$ from which it follows that $(\alpha_k x_k)_{k=0}^{\infty}\in L_p(\mathcal{M},\tau;c_0)$.

The proof for case (ii) is similar that of case (i), and so we omit it.
\end{proof}

\begin{cor}\label{NumberLimits}
Let $\alpha=(\alpha_k)_{k=0}^{\infty}\subset\mathbb{C}$ be such that there exists $L(\alpha)\in\mathbb{C}$ with $(\alpha_k-L(\alpha))_{k=0}^{\infty}\in c_0$ (i.e. $\alpha_k\to L(\alpha)$ in $\mathbb{C}$ as $k\to\infty$). If $(x_k)_{k=0}^{\infty}\in L_p(\mathcal{M},\tau;\ell_\infty)$ and $x\in L_p(\mathcal{M},\tau)$ are such that $(x_k-x)_{k=0}^{\infty}\in L_p(\mathcal{M},\tau;c_0)$, then $(\alpha_k x_k-L(\alpha)x)_{k=0}^{\infty}\in L_p(\mathcal{M},\tau;c_0)$.
\end{cor}

We will modify the notation for our averages now to include the weights. If $T\in DS^+(\mathcal{M},\tau)$, $\alpha=(\alpha_k)_{k=0}^{\infty}\subset\mathbb{C}$, and $n\in\mathbb{N}$, then write
$$M_n^{\alpha}(T):=\frac{1}{n}\sum_{k=0}^{n-1}\alpha_k T^k.$$

The next results show that the b.a.u. convergence of the weighted averages $(\frac{1}{n}\sum_{k=0}^{n-1}\alpha_k T^k(x_k))_{n=1}^{\infty}$ corresponds closely to the b.a.u. convergence of $(M_n^{\alpha}(T)(x))_{n=1}^{\infty}$ whenever $\alpha=(\alpha_k)_{k=0}^{\infty}\in\ell_\infty$ and $(x_k-x)_{k=0}^{\infty}\in L_p(\mathcal{M},\tau;c_0)$. Remarkably, almost no extra work is needed to prove this, as it only uses the already established results about the unweighted averages $(M_n(T))_{n\in\mathbb{N}}$ and the construction of more elements of $L_p(\mathcal{M},\tau;c_0)$ above. This time we will start with the $L_1$-case.

\begin{teo}\label{WeightedDirect}
Let $(\mathcal{M},\tau)$ be a tracial von Neumann algebra, $T:L_1(\mathcal{M},\tau)\to L_1(\mathcal{M},\tau)$ be induced by an ergodic normal $\tau$-preserving $*$-homomorphism of $\mathcal{M}$, and $\alpha=(\alpha_k)_{k=0}^{\infty}\in \ell_\infty$. If $(x_k)_{k=0}^{\infty}\in L_1(\mathcal{M},\tau;\ell_\infty)$ and $x\in L_1(\mathcal{M},\tau)$ are such that $(x_k-x)_{k=0}^{\infty}\in L_1(\mathcal{M},\tau;c_0)$, then $$\frac{1}{n}\sum_{k=0}^{n-1}\alpha_k T^k(x_k)-M_n^{\alpha}(T)(x)\to0\text{ b.a.u. as }n\to\infty.$$ Furthermore, if $M_n^{\alpha}(T)(x)\to F_T^{\alpha}(x)$ b.a.u. as $n\to\infty$ for some $F_T^{\alpha}(x)\in L_1(\mathcal{M},\tau)$, then $\frac{1}{n}\sum_{k=0}^{n-1}\alpha_k T^k(x_k)\to F_T^{\alpha}(x)$ b.a.u. as $n\to\infty$ as well.
\end{teo}
\begin{proof}
We observe that the linearity of each $T^k$ implies that
$$\frac{1}{n}\sum_{k=0}^{n-1}\alpha_k T^k(x_k)-M_n^{\alpha}(T)(x)=\frac{1}{n}\sum_{k=0}^{n-1}T^k\left(\alpha_k(x_k-x)\right)\text{ for every }n\in\mathbb{N}.$$ Since $(x_k-x)_{k=0}^{\infty}\in L_p(\mathcal{M},\tau;c_0)$ by assumption, it follows by Proposition \ref{BoundedWeightExample}(i) that $(\alpha_k (x_k-x))_{k=0}^{\infty}\in L_p(\mathcal{M},\tau;c_0)$ since $\alpha\in\ell_\infty$. Hence, the result follows immediately from Theorem \ref{DirectBreimanErgodicTheorem}. The statement about b.a.u. convergence follows from the linearity of b.a.u. convergence.
\end{proof}

The same proof could be applied verbatim to generalize Theorem \ref{GeneralBreimanErgodicTheorem} as well. However, we can show even more by instead proving that the family $(M_n^{\alpha}(T))_{n\in\mathbb{N}}$ is of strong type $(p,p)$ whenever $1<p<\infty$ and $\alpha$ is a (potentially unbounded) sequence in the set $W_s$ for an appropriate value of $s$, which we will now define.

\begin{comment}
\begin{teo}
Let $(\mathcal{M},\tau)$ be a semifinite von Neumann algebra, $1<p<\infty$, and $T\in DS^+(\mathcal{M},\tau)$.
Let $(x_k)_{k=0}^{\infty}\in L_p(\mathcal{M},\tau;\ell_\infty)$ and $x\in L_p(\mathcal{M},\tau)$ be such that $(x_k-x)_{k=0}^{\infty}\in L_p(\mathcal{M},\tau;c_0)$. If $\alpha=(\alpha_k)_{k=0}^{\infty}\in\ell_\infty$, then
$$\left(\frac{1}{n}\sum_{k=0}^{n-1}\alpha_k T^k(x_k)-M_n^{\alpha}(T)(x)\right)_{n=1}^{\infty}\in L_p(\mathcal{M},\tau;c_0).$$
Furthermore, if $M_n^{\alpha}(T)(x)\to F_T^{\alpha}(x)$ b.a.u. as $n\to\infty$ for some $F_T^{\alpha}(x)\in L_p(\mathcal{M},\tau)$, then $\frac{1}{n}\sum_{k=0}^{n-1}\alpha_k T^k(x_k)\to F_T^{\alpha}(x)$ b.a.u. as $n\to\infty$ as well. If the former b.a.u. convergence occurs as a consequence of the sequence of their differences being in $L_p(\mathcal{M},\tau;c_0)$, then the latter convergence occurs for the same reason.
\end{teo}
\end{comment}

Given $\alpha=(\alpha_k)_{k=0}^{\infty}\subset\mathbb{C}$ and $1\leq q<\infty$, define
$$\|\alpha\|_{W_q}:=\left(\limsup_{n\in\mathbb{N}}\frac{1}{n}\sum_{k=0}^{n-1}|\alpha_k|^q\right)^{1/q}$$ and let $$W_q:=\left\{\alpha=(\alpha_k)_{k=0}^{\infty}\subset\mathbb{C}:\|\alpha\|_{W_q}<\infty\right\}.$$ If $q=\infty$ then let $W_\infty=\ell_\infty$ with $\|\cdot\|_{W_\infty}=\|\cdot\|_{\ell_\infty}$. Then it is known that $W_q$ is a vector space over $\mathbb{C}$ with $\|\cdot\|_{W_q}$ as a seminorm. If $1\leq q<\infty$, then let $|\alpha|_{W_q}$ denote the same expression as $\|\alpha\|_{W_q}$ except with $\sup$ replacing $\limsup$, noting that $\alpha\in W_q$ if and only if $|\alpha|_{W_q}<\infty$. If $q=\infty$ then write $|\cdot|_{W_\infty}=\|\cdot\|_{W_\infty}$. If $1\leq q<s\leq\infty$ then $\|\alpha\|_{W_q}\leq\|\alpha\|_{W_s}$ and $|\alpha|_{W_q}\leq|\alpha|_{W_s}$, either of which implies that $W_s\subset W_q$

It was shown in \cite{ob4} that, if $T\in DS^+(\mathcal{M},\tau)$ and $1<p,q<\infty$ are such that $\frac{1}{p}+\frac{1}{q}=1$, then $(|\alpha|_{W_q}^{-1}M_n^{\alpha}(T))_{(n,\alpha)\in\mathbb{N}\times W_q}$ satisfies a weak type $(p,p)$ maximal inequality (see \cite{ob4} for more details on a definition of this).  It is known that if a family is of strong type $(p,p)$, then it is of weak type $(p,p)$ as well; the converse does not hold in general. Unfortunately, the averages $(M_n^{\alpha}(T))_{n\in\mathbb{N}}$ are not necessarily of strong type $(p,p)$ in general when $\alpha\in W_q$. However, if $\alpha$ is in a slightly smaller $W_s$ space then the strong type maximal inequality will hold on $L_p(\mathcal{M},\tau)$, as we will now show.

\begin{teo}\label{WeightedMaximalInequality}
Let $(\mathcal{M},\tau)$ be a semifinite von Neumann algebra, $T\in DS^+(\mathcal{M},\tau)$, and $1<p,q<\infty$ be such that $\frac{1}{p}+\frac{1}{q}=1$. If $s>q$, $\alpha\in W_s$, and $1<r<\infty$ is such that $\frac{1}{r}+\frac{1}{s}=1$, then the family $(M_n^{\alpha}(T))_{n\in\mathbb{N}}$ is of strong type $(p,p)$ with constant $16C_{p/r}^{1/r}|\alpha|_{W_s}$, where $C_{p/r}$ is the strong type constant of $(M_n(T))_{n\in\mathbb{N}}$ acting on $L_{p/r}(\mathcal{M},\tau)$.
\end{teo}
\begin{proof}
We first assume that $\alpha=(\alpha_k)_{k=0}^{\infty}\in W_s$ is such that $\alpha_k\geq0$ for every $k\in\mathbb{N}_0$.
%Since $\alpha=(\alpha_k)_{k=0}^{\infty}\in W_s$ if and only if $\alpha^*:=(\overline{\alpha_k})_{k=0}^{\infty}\in W_s$ and $|\alpha|_{W_s}=|\alpha^*|_{W_s}$, writing $Re(\alpha):=\frac{1}{2}(\alpha+\alpha^*),Im(\alpha):=\frac{1}{2}(\alpha-\alpha)^*\in W_s$, without loss of generality we may assume that $\alpha_k\in\mathbb{R}$ for every $k\in\mathbb{N}_0$ by instead proving the inequality for $Re(\alpha)$ and $Im(\alpha)$ and using linearity since $Re(\alpha),Im(\alpha)\subset\mathbb{R}$,
%$$(M_n^\alpha(T)(x))_{n=1}^{\infty}=(M_n^{Re(\alpha)}(T)(x))_{n=1}^{\infty}+i(M_n^{Im(\alpha)}(T)(x))_{n=1}^{\infty}$$ for every $x\in L_p(\mathcal{M},\tau)$, and
%$|Re(\alpha)|_{W_s},|Im(\alpha)|_{W_s}\leq|\alpha|_{W_s}.$ Similarly, by writing $\alpha_k=\beta_k-\gamma_k$ with $\beta_k,\gamma_k\in[0,\infty)$, if we denote $\beta:=(\beta_k)_{k=0}^{\infty}$ and $\gamma:=(\gamma_k)_{k=0}^{\infty}$, then $|\beta|_{W_s},|\gamma|_{W_s}\leq|\alpha|_{W_s}$, so that $\beta,\gamma\in W_s$ and $\alpha=\beta-\gamma$. By proving the result first for $\beta$ and $\gamma$ and then use the fact that $$(M_n^{\alpha}(T)(x))_{n=1}^{\infty}=(M_n^{\beta}(T)(x))_{n=1}^{\infty}-(M_n^{\gamma}(T)(x))_{n=1}^{\infty}$$ and $|\beta|_{W_s},|\gamma|_{W_s}\leq|\alpha|_{W_s}$, without loss of generality we may assume that $\alpha_k\geq0$ for every $k\in\mathbb{N}_0$. Note then that $\alpha_k T^k$ is a positive linear operator of $L_t(\mathcal{M},\tau)$ for every $k\in\mathbb{N}_0$ and $1\leq t\leq\infty$ (and thus so is $M_n^{\alpha}(T)$ for every $n\in\mathbb{N}$).
Fix $x\in L_1(\mathcal{M},\tau)\cap\mathcal{M}^{+}$. It follows by \cite[Lemma 3.1]{ob4} that
$$\frac{1}{n}\sum_{k=0}^{n-1}\alpha_k T^k(x)\leq\left(\frac{1}{n}\sum_{k=0}^{n-1}\alpha_k^s\right)^{1/s}\left(\frac{1}{n}\sum_{k=0}^{n-1}T^k(x^r)\right)^{1/r}\leq|\alpha|_{W_s}M_n^{\alpha}(T)(x^r)^{1/r}.$$
Since $x\in L_p(\mathcal{M},\tau)^{+}$ it follows that $x^r\in L_{p/r}(\mathcal{M},\tau)$. Since $p/r>1$ it follows that $(M_n(T))_{n\in\mathbb{N}}$ is of strong type $(p/r,p/r)$ with constant $C_{p/r}>0$; hence there exists $c\in L_{p/r}(\mathcal{M},\tau)^{+}$ with
$$\|c\|_{p/r}\leq C_{p/r}\|x^r\|_{p/r}\text{ and }M_n^{\alpha}(T)(x^r)\leq c\text{ for every }n\in\mathbb{N}.$$
If $a\in L_{p/r}^{+}$ one can find that $\|a\|_{p/r}=\|a^{1/r}\|_{p}^{r}$, so that the first inequality shows that $\|c^{1/r}\|_{p}\leq C_{p/r}^{1/r}\|x\|_p$. The second inequality implies that
$$|\alpha|_{W_s}M_n^{\alpha}(T)(x^r)^{1/r}\leq|\alpha|_{W_s}c^{1/r}$$ since the map $0\leq t\mapsto t^{1/r}$ is increasing with respect to the ordering on $L_0(\mathcal{M},\tau)^{+}$ by \cite[Theorem 5]{dt}. Hence
$$0\leq\frac{1}{n}\sum_{k=0}^{n-1}\alpha_k T^k(x)\leq|\alpha|_{W_s}c^{1/r}\text{ for every }n\in\mathbb{N}$$ with $\||\alpha|_{W_s}c^{1/r}\|_p\leq |\alpha|_{W_s}C_{p/r}^{1/r}\|x\|_p.$ It follows by the characterization of positive sequences in $L_p(\mathcal{M},\tau;\ell_\infty)$ that $(M_n^{\alpha}(T)(x))_{n=1}^{\infty}\in L_p(\mathcal{M},\tau;\ell_\infty)$ with a norm of at most $C_{p/r}^{1/r}|\alpha|_{W_s}\|x\|_p$.

By writing $x=\sum_{j=0}^{3}i^j x_j$, where $x_j\in L_1(\mathcal{M},\tau)\cap\mathcal{M}^{+}$ and $\|x_j\|_p\leq\|x\|_p$ for each $j\in\{0,1,2,3\}$, and applying the triangle inequality, it follows that 
$$
(M_n^{\alpha}(T)(x))_{n=1}^{\infty}
=\sum_{j=0}^{3}i^j(M_n^{\alpha}(T)(x_j))_{n=1}^{\infty}\in L_p(\mathcal{M},\tau;\ell_\infty)
$$ 
and 
\begin{align*}
\left\|{\sup_{k\in\mathbb{N}_0}}^{+}M_n^{\alpha}(T)(x)\right\|_p
\leq&
\sum_{j=0}^{3}\left\|{\sup_{k\in\mathbb{N}_0}}^{+}M_n^{\alpha}(T)(x_j)\right\|_p \\
\leq&\sum_{j=0}^{3}|\alpha|_{W_s}C_{p/r}^{1/r}\|x_j\|_p
\leq4|\alpha|_{W_s}C_{p/r}^{1/r}\|x\|_p.
\end{align*}

Since $x\in L_1(\mathcal{M},\tau)\cap\mathcal{M}$ was arbitary, it follows that $(M_n^{\alpha}(T))_{n\in\mathbb{N}}$ defines a bounded linear map from the dense subspace $L_1(\mathcal{M},\tau)\cap\mathcal{M}$ of $L_p(\mathcal{M},\tau)$ to $L_p(\mathcal{M},\tau;\ell_\infty)$. Hence, its unique continuous linear extension is bounded with the same norm and can be seen to agree with taking $x\in L_p(\mathcal{M},\tau)$ to $(M_n^{\alpha}(T)(x))_{n=1}^{\infty}$, meaning $(M_n^{\alpha}(T))_{n\in\mathbb{N}}$ is of strong type $(p,p)$ with constant $4C_{p/r}^{1/r}|\alpha|_{W_s}$ when $\alpha=(\alpha_k)_{k=0}^{\infty}$ is such that $\alpha_k\geq0$ for every $k\in\mathbb{N}_0$.
 
Similar to the above, if $\alpha\in W_s$ is arbitrary then by taking the real and imaginary parts of $\alpha$ and splitting those up in to their positive and negative parts one may write $\alpha=\sum_{j=0}^{3}i^j \alpha^{(j)}$, where $\alpha^{(j)}=(\alpha_k^{(j)})_{k=0}^{\infty}\in W_s$ is such that $\alpha_k^{(j)}\geq0$ for every $k\in\mathbb{N}_0$ and $|\alpha^{(j)}|_{W_s}\leq|\alpha|_{W_s}$ for each $j\in\{0,1,2,3\}$. If $x\in L_p(\mathcal{M},\tau)$, then
$$
\left(M_n^{\alpha}(T)(x)\right)_{n=1}^{\infty}
=\left(\sum_{j=0}^{3}M_n^{\alpha^{(j)}}(T)(x)\right)_{n=1}^{\infty}
\in L_p(\mathcal{M},\tau;\ell_\infty)$$ and
\begin{align*}
\left\|{\sup_{k\in\mathbb{N}_0}}^{+}M_n^{\alpha}(T)(x)\right\|_p
\leq&\sum_{j=0}^{3}\left\|{\sup_{k\in\mathbb{N}_0}}^{+}M_n^{\alpha^{(j)}}(T)(x)\right\|_p \\
\leq&\sum_{j=0}^{3}4|\alpha^{(j)}|_{W_s}C_{p/r}^{1/r}\|x\|_p
\leq16|\alpha|_{W_s}C_{p/r}^{1/r}\|x\|_p
\end{align*}
finishing the proof for arbitrary $\alpha\in W_s$.
\end{proof}

One can then repeat the proof of Theorem \ref{GeneralBreimanErgodicTheorem} by replacing $(M_n(T))_{n\in\mathbb{N}}$ with $(M_n^{\alpha}(T))_{n\in\mathbb{N}}$ and first considering $\alpha=(\alpha_k)_{k=0}^{\infty}\in W_s$ with $\alpha_k\geq0$ for every $k\in\mathbb{N}_0$, and then taking linear combinations of the separate parts of $\alpha$ (as was done in the previous proof) to obtain the following weighted version of Theorem \ref{WeightedDirect}.

\begin{teo}\label{WeightedGeneralized}
Let $(\mathcal{M},\tau)$ be a semifinite von Neumann algebra, $T\in DS^+(\mathcal{M},\tau)$, and $1<p,q<\infty$ be such that $\frac{1}{p}+\frac{1}{q}=1$. Let $s>q$ and $\alpha\in W_s$. Then if $(x_k)_{k=0}^{\infty}\in L_p(\mathcal{M},\tau;\ell_\infty)$ and $x\in L_p(\mathcal{M},\tau)$ are such that $(x_k-x)_{k=0}^{\infty}\in L_p(\mathcal{M},\tau;c_0)$, then
$$\left(\frac{1}{n}\sum_{k=0}^{n-1}\alpha_k T^k(x_k)-M_n^{\alpha}(T)(x)\right)_{n=1}^{\infty}\in L_p(\mathcal{M},\tau;c_0).$$
Furthermore, if $M_n^{\alpha}(T)(x)\to F_T^{\alpha}(x)$ b.a.u. as $n\to\infty$ for some $F_T^{\alpha}(x)\in L_p(\mathcal{M},\tau)$, then $\frac{1}{n}\sum_{k=0}^{n-1}\alpha_k T^k(x)\to F_T^{\alpha}(x)$ b.a.u. as $n\to\infty$ as well. 
\end{teo}

We will now quickly describe a few classes of weights for which these results are known to apply to.

\begin{ex}
A commonly studied family of weights in ergodic theory are the Besicovitch sequences. A sequence $P=(P_k)_{k=0}^{\infty}\subset\mathbb{C}$ is called a \textit{trigonometric polynomial} if there exists $c_1,...,c_r\in\mathbb{C}$ and $\lambda_1,...,\lambda_r\in\mathbb{T}$ such that
$$P_k=\sum_{s=1}^{r}c_s\lambda_s^k$$ for every $k\in\mathbb{N}_0$. For a fixed $1<q<\infty$, a sequence $\alpha=(\alpha_k)_{k=0}^{\infty}\subset\mathbb{C}$ is called a \textit{$q$-Besicovitch sequence} if for every $\epsilon>0$ there exists a trigonometric polynomial $P_{\epsilon}=(P_{\epsilon,k})_{k=0}^{\infty}$ such that
$$\|\alpha-P_\epsilon\|_{W_q}<\epsilon.$$
Write $B_q$ for the set of all such sequences. Write $B_\infty:=B_1\cap W_\infty$, calling elements of this set \textit{bounded Besicovitch sequences}. One can find, similar to $W_q$, that $B_s\subset B_q$ whenever $1\leq q<s\leq\infty$.
	
If $(\mathcal{M},\tau)$ is a semifinite von Neumann algebra with a separable predual, $T\in DS^+(\mathcal{M},\tau)$, either $1<p,q<\infty$ are such that $\frac{1}{p}+\frac{1}{q}=1$ or $p=1$ and $q=\infty$, and $\alpha\in B_q$, then it was shown that $(M_n^{\alpha}(T)(x))_{n=1}^{\infty}$ converges b.a.u. as $n\to\infty$  for every $x\in L_p(\mathcal{M},\tau)$ for the $p=1$ case in \cite[Theorem 4.6]{cls} and for the $1<p<\infty$ case in \cite{ob4}. In the $1<p<\infty$ case Theorem \ref{WeightedGeneralized} would apply whenever $\alpha\in B_s$ for some $s>q$ since then $B_s\subseteq B_q$.
\end{ex}

\begin{ex}
Let $(\Omega,\mathcal{F},\mu)$ be a probability space, and let $\mathcal{B}_{\mathbb{C}}$ denote the Borel $\sigma$-algebra of $\mathbb{C}$. A family $(Y_k)_{k=0}^{\infty}$ of measurable functions $Y_k:\Omega\to\mathbb{C}$ is called an i.i.d. sequence on $\Omega$ if it satisfies:
\begin{itemize}
\item[(i)](Independent) $\mu\left(\bigcap_{j=1}^{m}Y_{k_j}^{-1}(E_j)\right)=\prod_{k=1}^{m}\mu(Y_{k_j}^{-1}(E_j))$ for every $E_1,...,E_m\in\mathcal{B}_{\mathbb{C}}$, every $k_1,...,k_m\in\mathbb{N}_0$, and every $m\in\mathbb{N}$,
\item[(ii)](Identically Distributed) $\mu\circ Y_k^{-1}=\mu\circ Y_j^{-1}$ as measures on $(\mathbb{C},\mathcal{B}_{\mathbb{C}})$ for every $k,j\in\mathbb{N}_0$.
\end{itemize}
	
It was shown in \cite[Theorems 3.2 and 3.3]{ob5} that if either $p=1$ and $q=\infty$ or if $1<p,q<\infty$ are such that $\frac{1}{p}+\frac{1}{q}=1$, and if $Y_0\in L_q(\Omega,\mathcal{F},\mu)$, then there exists $\Omega^{\prime}\subseteq\Omega$ with $\mu(\Omega^{\prime})=1$ such that if $\omega\in\Omega^{\prime}$ and $Y(\omega):=(Y_k(\omega))_{k=0}^{\infty}$, then $Y(\omega)\in W_q$ and $$\frac{1}{n}\sum_{k=0}^{n-1}Y_k(\omega)T^k(x)\to F_T(x)\int_{\Omega}Y_0\ d\mu\text{ b.a.u. as }n\to\infty$$ for every semifinite von Neumann algebra $(\mathcal{M},\tau)$, $T\in DS^+(\mathcal{M},\tau)$, and $x\in L_p(\mathcal{M},\tau)$.

If $1<q<s\leq\infty$, then $(\Omega,\mathcal{F},\mu)$ being a probability space implies that $L_s(\Omega,\mathcal{F},\mu)\subset L_q(\Omega,\mathcal{F},\mu)$, so that the above holds on $L_p(\mathcal{M},\tau)$ for $Y(\omega)$ as defined above if $Y_0\in L_s(\Omega,\mathcal{F},\mu)$ since this implies $\|Y_0\|_q\leq\|Y_0\|_s<\infty$ and $Y(s)\in W_s\subset W_q$.
\end{ex}

\begin{ex}
Let $\mu:\mathbb{N}_0\to\mathbb{C}$ denote the number-theoretic M\"{o}bius function and $\lambda:\mathbb{N}_0\to\mathbb{C}$ denote the Liouville function. Then \cite[Theorem 3.4]{ob5} implies that if $(\mathcal{M},\tau)$ is a semifinite von Neumann algebra, $1\leq p<\infty$, and $T\in DS^+(\mathcal{M},\tau)$, then both
$$\frac{1}{n}\sum_{k=0}^{n-1}\mu(k)T^k(x)\text{ and }\frac{1}{n}\sum_{k=0}^{n-1}\lambda(k)T^k(x)$$ converge to $0$ b.a.u. as $n\to\infty$ for every $x\in L_p(\mathcal{M},\tau)$. (Note that both $(\mu(k))_{k=0}^{\infty}$ and $(\lambda(k))_{k=0}^{\infty}$ are in $W_\infty$ due to the ranges of $\mu$ and $\lambda$ being subsets of $\{-1,0,1\}$.) A generalization to this will be discussed below in Theorem \ref{MobiusPolynomialErgodicTheorem}.
\end{ex}

\subsection{Subsequential Averages}% and other generalizations}

%\begin{rem}\label{MoreGeneralSetting}
Many of the results in this article can be generalized further than stated. Namely, the main properties of $T\in DS^+(\mathcal{M},\tau)$ that were used in this article so far were:
\begin{itemize}
\item[(i)] $T$ is positive,
\item[(ii)] $(M_n(T))_{n\in\mathbb{N}}$ is of strong type $(p,p)$,
\item[(iii)] for every $x\in L_p(\mathcal{M},\tau)$ there exists $F_T(x)\in L_p(\mathcal{M},\tau)$ such that $$(M_n(T)-F_T(x))_{n=1}^{\infty}\in L_p(\mathcal{M},\tau;c_0).$$
\end{itemize}
The proofs of Theorems \ref{GeneralBreimanErgodicTheorem} and \ref{WeightedGeneralized} would follow immediately if $T\in DS^+(\mathcal{M},\tau)$ was replaced a positive bounded linear operator $T:L_p(\mathcal{M},\tau)\to L_p(\mathcal{M},\tau)$ for a fixed $1\leq p<\infty$ satisfying (i) and (ii), with (iii) only being needed to get the b.a.u. convergence when $(x_k-x)_{k=0}^{\infty}\in L_p(\mathcal{M},\tau;\ell_\infty)$.  Since $(\frac{1}{n}\sum_{k=0}^{n-1}T^k(x_k))_{n=1}^{\infty}\in L_p(\mathcal{M},\tau;c_0)$ for every $(x_k)_{k=0}^{\infty}\in L_p(\mathcal{M},\tau;c_0)$ only needs conditions (i) and (ii) to prove, if b.a.u. convergence is the primary goal then (iii) in this list could also be replaced by the more general condition:
\begin{itemize}
\item[(iv)] for every $x\in L_p(\mathcal{M},\tau)$ there exists $F_T(x)\in L_p(\mathcal{M},\tau)$ such that $M_n(T)(x)\to F_T(x)$ b.a.u. as $n\to\infty$.
\end{itemize}
As a consequence, however, one would only have that $M_n(T)(x)\to F_T(x)$ b.a.u. as $n\to\infty$ if they only assumed (iv) and not necessarily $(M_n(T)(x)-F_T(x))_{n=1}^{\infty}\in L_p(\mathcal{M},\tau;c_0)$.
%\end{rem}

Another way of generalizing the averages $(M_n(T))_{n\in\mathbb{N}}$ is to replace them with the subsequential averages
$$M_{\textbf{k},n}(T):=\frac{1}{n}\sum_{j=0}^{n-1}T^{k_j} \text{ for every }n\in\mathbb{N},$$ where $\textbf{k}=(k_j)_{j=0}^{\infty}\subset\mathbb{N}_0$ is a (fixed) strictly increasing sequence.

The comments above also indicate that Theorems \ref{GeneralBreimanErgodicTheorem} and \ref{WeightedGeneralized} can be generalized to apply to these subsequential averages if one replaces $(M_n(T))_{n\in\mathbb{N}}$ by $(M_{\textbf{k},n}(T))_{n\in\mathbb{N}}$ in conditions (ii), (iii), and (iv). We summarize this in the next result. We omit the case $p=1$ since strong type $(1,1)$ maximal ergodic inequalities fairly often do not hold for subsequential averages. However, as was the case considered in Theorem \ref{DirectBreimanErgodicTheorem}, if it holds that $M_{\textbf{k},n}(T)(x)\to \tau(x)\textbf{1}$ b.a.u. as $n\to\infty$ for every $x\in L_1(\mathcal{M},\tau)$ then the corresponding version of Theorem \ref{DirectBreimanErgodicTheorem} would hold for $(M_{\textbf{k},n}(T))_{n\in\mathbb{N}}$ as well with the same proof.

\begin{teo}\label{SubsequenceResult}
Let $(\mathcal{M},\tau)$ be a semifinite von Neumann algebra, $1<p<\infty$, and $T\in DS^+(\mathcal{M},\tau)$. Let $\textbf{k}=(k_j)_{j=0}^{\infty}\subset\mathbb{N}_0$ be a strictly increasing sequence such that $(M_{\textbf{k},n}(T))_{n\in\mathbb{N}}$ is of strong type $(p,p)$. Then if $(x_k)_{k=0}^{\infty}\in L_p(\mathcal{M},\tau;c_0)$ it follows that
$$\left(\frac{1}{n}\sum_{j=0}^{n-1}T^{k_j}(x_j)\right)_{n=1}^{\infty},\left(\frac{1}{n}\sum_{j=0}^{n-1}T^{k_j}(x_{k_j})\right)_{n=1}^{\infty}\in L_p(\mathcal{M},\tau;c_0).$$

In addition to this the operator $(x_k)_{k=0}^{\infty}\mapsto\left(\frac{1}{n}\sum_{j=0}^{n-1}T^{k_j}(x_j)\right)_{n=1}^{\infty}$ from $L_p(\mathcal{M},\tau;\ell_\infty)$ to itself is bounded and leaves $L_p(\mathcal{M},\tau;c_0)$ invariant.

Furthermore, if $(x_k)_{k=0}^{\infty}\in L_p(\mathcal{M},\tau;\ell_\infty)$ and $x\in L_p(\mathcal{M},\tau)$ are such that $(x_k-x)_{k=0}^{\infty}\in L_p(\mathcal{M},\tau;c_0)$, and if $M_{\textbf{k},n}(T)(x)\to F_{T,\textbf{k}}(x)$ b.a.u. as $n\to\infty$ for some $F_{T,\textbf{k}}\in L_p(\mathcal{M},\tau)$, then it follows that $$\frac{1}{n}\sum_{j=0}^{n-1}T^{k_j}(x_j)\to F_{T,\textbf{k}}(x)\text{ and } \frac{1}{n}\sum_{j=0}^{n-1}T^{k_j}(x_{k_j})\to F_{T,\textbf{k}}(x)$$ b.a.u. as $n\to\infty$ as well. If the b.a.u. convergence of $M_{\textbf{k},n}(T)(x)\to F_{T,\textbf{k}}(x)$ occurs due to the sequence of differences being in $L_p(\mathcal{M},\tau;c_0)$, then the latter convergences occur for the same reason. 
\end{teo}

In addition to the discussion above, the boundedness of the map $(x_k)_{k=0}^{\infty}\mapsto\left(\frac{1}{n}\sum_{j=0}^{n-1}T^{k_j}(x_j)\right)_{n=1}^{\infty}$ follows the same argument of Corollary \ref{Corollary}.

We will now discuss a few subsequences of $\mathbb{N}$ which are good for the individual ergodic theorem on $L_p(\mathcal{M},\tau)$ in the sense that $(M_{\textbf{k},n}(T)(x))_{n=1}^{\infty}$ converges b.a.u. for every $x\in L_p(\mathcal{M},\tau)$ for every $T\in DS^+(\mathcal{M},\tau)$ (or at least for $T$ induced by a normal $\tau$-preserving $*$-automorphism of $\mathcal{M}$).

If $\textbf{k}=(k_j)_{j=0}^{\infty}$ is an increasing sequence in $\mathbb{N}_0$ then define its natural density $d(\textbf{k})$ to be the quanitity
$$d(\textbf{k}):=\lim_{n\to\infty}\frac{\text{card}\left(\{k_j:j\in\mathbb{N}_0\}\cap[0,n]\right)}{n+1},$$ assuming that the limit exists. 

\begin{ex}
	If $(\mathcal{M},\tau)$ is a semifinite von Neumann algebra, $1\leq p<\infty$, $T\in DS^+(\mathcal{M},\tau)$, and $d(\textbf{k})=1$, then $M_{\textbf{k},n}(T)(x)\to F_T(x)$ b.a.u. as $n\to\infty$ for every $x\in L_p(\mathcal{M},\tau)$ by \cite[Theorem 3.1]{ob1} and $(M_{\textbf{k},n}(T))_{n\in\mathbb{N}}$ can be seen to be of strong type $(p,p)$ if $1<p<\infty$ by the proofs of \cite[Proposition 3.1]{ob1} and \cite[Theorem 3.1]{mmt} since these show that the subsequential averages $(M_{\textbf{k},n}(T))_{n\in\mathbb{N}}$ can be converted to weighted averages by bounded weights which are strong type $(p,p)$ by Theorem \ref{WeightedMaximalInequality} and the fact that $W_\infty\subset W_s$ for every $1\leq s<\infty$.
\end{ex}

\begin{ex}
	Let $(Y_k)_{k=0}^{\infty}$ be a sequence of i.i.d. random variables on a probability space $(\Omega,\mathcal{F},\mu)$ such that $Y_k(\omega)\in\{0,1\}$ for every $k\in\mathbb{N}_0$ and $\omega\in\Omega$, and assume that $\mu(Y_0^{-1}(\{1\}))>0$. For any $\omega\in\Omega$ define $k_0(\omega)=\inf\{k\in\mathbb{N}_0:Y_k(\omega)=1\}$ and
	$$k_j(\omega)=\inf\{k\in\mathbb{N}:k\geq k_{j-1}\text{ and }Y_k(\omega)=1\}\text{ for }j\in\mathbb{N}.$$
	Then for $\mu$-a.e. $\omega\in\Omega$ the sequence $\textbf{k}(\omega):=(k_j)_{j=0}^{\infty}$ has natural density $d(\textbf{k}(\omega))=\mu(Y_0^{-1}(\{1\}))$, and so it follows from \cite[Theorem 3.2]{ob5} and \cite[Proposition 3.2]{ob1} that if $(\mathcal{M},\tau)$ is a semifinite von Neumann algebra, $1\leq p<\infty$, and $T\in DS^+(\mathcal{M},\tau)$ then $M_{\textbf{k}(\omega),n}(T)(x)\to F_T(x)$ b.a.u. as $n\to\infty$, and $(M_{\textbf{k}(\omega),n}(T))_{n\in\mathbb{N}}$ is of strong type $(p,p)$ if $1<p<\infty$.
\end{ex}

\begin{ex}
	An variant of the previous example is known as well. Let $0<\alpha<1$, $(\Omega,\mathcal{F},\mu)$ be a probability space, and let $(Y_k)_{k=1}^{\infty}$ be an independent family of random variables from $\Omega$ to $\{0,1\}$ such that $\mu(Y_n^{-1}(\{1\}))=n^{-\alpha}$ for every $n\in\mathbb{N}$.
	
	Define functions $n_k:\Omega\to\mathbb{N}$ for $k\geq1$ defined by
	$$n_k(\omega):=\min\{N\geq1:\sum_{j=1}^{N}Y_k(\omega)=k\}.$$
	One can find that $d((n_k(\omega))_{k=0}^{\infty})=0$ for $\mu$-a.e. $\omega\in\Omega$.
	
	Let $(\mathcal{M},\tau)$ be a semifinite von Neumann algebra, $1<p<\infty$, and $T\in DS^+(\mathcal{M},\tau)$. It was shown in \cite[Theorem 5.2]{lmz} that there exists $\Omega^{\prime}\in\mathcal{F}$ with $\mu(\Omega^{\prime})=1$ such that $\omega\in\Omega^{\prime}$ implies that $$\frac{1}{m}\sum_{k=1}^{m}T^{n_k(\omega)}(x)\to F_T(x)\text{ b.a.u. as }n\to\infty$$ for every $x\in L_p(\mathcal{M},\tau)$.
\end{ex}

\begin{ex}
	Let $(\mathcal{M},\tau)$ be a semifinite von Neumann algebra and $T:\mathcal{M}\to\mathcal{M}$ be a normal $\tau$-preserving $*$-homomorphism (and use $T$ for its extensions to $L_p(\mathcal{M},\tau)$ for every $1\leq p<\infty$). Fix $1<p<\infty$.
	
	%If $\textbf{k}^{t}:=(k^t)_{k=0}^{\infty}$ for some fixed $t\in\mathbb{N}$ with $t\geq2$, then it was shown in \cite[Theorem 1.1]{chw} that $(M_{\textbf{k}^{t},n}(T))_{n\in\mathbb{N}}$ is of strong type $(p,p)$ whenever $p>\frac{1+\sqrt{5}}{2}$. It was recently shown in \cite[Remark 4.9]{hw} that the averages converge b.a.u. as $n\to\infty$ for every $x\in L_p(\mathcal{M},\tau)$ for the same values of $p$. 
	
	It was shown in \cite{hlrw} that, if $P\in\mathbb{Z}[X]$ is a polynomial with integer coefficients which maps $\mathbb{N}_0$ to itself, then $P(\textbf{k}):=(P(k))_{k=0}^{\infty}$ is of strong type $(p,p)$ and for every $x\in L_p(\mathcal{M},\tau)$ there exists $y\in L_p(\mathcal{M},\tau)$ such that $M_{P(\textbf{k}),n}(T)(x)\to y$ b.a.u. as $n\to\infty$.
	
	If $\textbf{p}=(p_k)_{k=0}^{\infty}$, where $p_k$ denotes the $k$-th prime number, then \cite{hw} shows that $(M_{\textbf{p},n}(T))_{n\in\mathbb{N}}$ is of strong type $(p,p)$ for every $1<p<\infty$ and that for every $x\in L_p(\mathcal{M},\tau)$ there exists $y\in L_p(\mathcal{M},\tau)$ such that $M_{\textbf{p},n}(T)(x)\to y$ b.a.u. as $n\to\infty$.
\end{ex}

Weighted versions of Theorem \ref{SubsequenceResult} can be proven as well, though the convergence of the weighted subsequential averages would still need to be proven. A number of such sequences were studied in \cite{ob1}. We note that a condition similar to a condition on the weights considered in \cite{ob5} holds for the subsequence $(k^t)_{k=0}^{\infty}$ when weighted by the number-theoretic M\"{o}bius function $(\mu(k))_{k=0}^{\infty}$, which we will now prove.

\begin{teo}\label{MobiusPolynomialErgodicTheorem}
	Let $(\mathcal{M},\tau)$ be a tracial von Neumann algebra, $1<p<\infty$, and $T:L_p(\mathcal{M},\tau)\to L_p(\mathcal{M},\tau)$ be induced by a normal $\tau$-preserving $*$-automorphism of $\mathcal{M}$. If $t\in\mathbb{N}$ is fixed and $\textbf{k}^t:=(k^t)_{k=0}^{\infty}$, then if $\mu:\mathbb{N}_0\to\{-1,0,1\}$ denotes the M\"{o}bius function then the weighted subsequential averages
	$$M_{\textbf{k}^t,n}^{\mu}(T)(x):=\frac{1}{n}\sum_{k=0}^{n-1}\mu(k)T^{k^t}(x)$$ converge b.a.u. as $n\to\infty$ for every $x\in L_p(\mathcal{M},\tau)$. The b.a.u. limit is $0$ for all $x\in L_p(\mathcal{M},\tau)$ if either $p=2$ or if $1<p<\infty$ and $\tau(\textbf{1})<\infty$.	
	
	Furthermore, if $(x_k)_{k=0}^{\infty}\in L_p(\mathcal{M},\tau;\ell_\infty)$ and $x\in L_p(\mathcal{M},\tau)$ are such that $(x_k-x)_{k=0}^{\infty}\in L_p(\mathcal{M},\tau;c_0)$, then
	$$\frac{1}{n}\sum_{k=0}^{n-1}\mu(k)T^{k^t}(x_k)$$ converge b.a.u. to the same limit as $(M_{\textbf{k}^t,n}^{\mu}(T)(x))_{n=1}^{\infty}$ as $n\to\infty$.
\end{teo}
\begin{proof}
	It was stated in \cite[Inequality (1.3)]{jl} that the methods of \cite[Theorem 10]{hua} imply that for every $A>0$ and $t\in\mathbb{N}$ there exists $C_{A,t}>0$ such that
	$$\sup_{\theta\in\mathbb{R}}\left|\sum_{k\leq x}\mu(k)e^{2\pi i k^{t}\theta}\right|\leq\frac{C_{A,t}x}{\ln(x)^{A}}$$
	for every $x>0$. From this it can be shown that
	$$\sum_{m=1}^{\infty}\sup_{\lambda\in\mathbb{T}}\left|\frac{1}{\rho^m}\sum_{k=0}^{[\rho^m]-1}\mu(k)\lambda^{k^t}\right|^2
	\leq\frac{C_{1,t}^2}{(\ln(\rho))^{2}}\sum_{m=1}^{\infty}\frac{1}{m^2}<\infty$$ 
	for every $\rho>1$. Since the averages $(M_{\textbf{k}^t,n}(T)(x))_{n=1}^{\infty}$ converge b.a.u. as $n\to\infty$ for every $x\in L_p(\mathcal{M},\tau)$ by \cite[Theorem 1.1]{hlrw}, the b.a.u. convergence follows by repeating the proof of \cite[Theorem 3.1]{ob5} with every instance of $M_n(T)$ replaced by $M_{\textbf{k}^t,n}(T)$. 
	
	A weak type $(p,p)$ maximal inequality for the averages $(M_{\textbf{k}^{t},n}^{\mu}(T))_{n\in\mathbb{N}}$ is needed in \cite[Theorem 3.1]{ob5} to extend the convergence from $L_2(\mathcal{M},\tau)$ to $L_p(\mathcal{M},\tau)$, but that would be implied by a corresponding strong type $(p,p)$ maximal inequality. The necessary strong type maximal inequality is a consequence of the continuity of the composition of the maps $$x\mapsto(x)_{k\in\mathbb{N}_0}\mapsto(\mu(k)x)_{k\in\mathbb{N}_0}\mapsto\left(\frac{1}{n}\sum_{k=0}^{n-1}T^{k^t}(\mu(k)x)\right)_{n=1}^{\infty}=\left(\frac{1}{n}\sum_{k=0}^{n-1}\mu(k)T^{k^t}(x)\right)_{n=1}^{\infty},$$ where the first map takes $L_p(\mathcal{M},\tau)$ to $L_p(\mathcal{M},\tau;\ell_\infty)$ and the other maps all take the latter space to itself. The second map is continuous by Proposition \ref{BoundedWeightExample}, while the third map is continuous by Theorem \ref{SubsequenceResult}. The first map from $L_p(\mathcal{M},\tau)$ to $L_p(\mathcal{M},\tau;\ell_\infty)$ is continuous since writing $x\in L_p(\mathcal{M},\tau)$ with $x=u|x|$ as its polar decomposition implies that $x=(u|x|^{1/2})\textbf{1}(|x|^{1/2})$ with $|x|^{1/2},u|x|^{1/2}\in L_{2p}(\mathcal{M},\tau)$, which means that $$\left\|{\sup_{k\in\mathbb{N}_0}}^{+}x\right\|_p\leq\|u|x|^{1/2}\|_{2p}\|\textbf{1}\|_\infty\||x|^{1/2}\|_{2p}\leq\||x|^{1/2}\|_{2p}^{2}=\|x\|_p$$ for every $x\in L_p(\mathcal{M},\tau)$ by the infimum definition of the norm of $L_p(\mathcal{M},\tau;\ell_\infty)$.
	
	If $(x_k)_{k=0}^{\infty}\in L_p(\mathcal{M},\tau;\ell_\infty)$ is such that $(x_k-x)_{k=0}^{\infty}\in L_p(\mathcal{M},\tau;c_0)$, then it follows that $(\mu(k)(x_k-x))_{k=0}^{\infty}\in L_p(\mathcal{M},\tau;c_0)$ as well by Proposition \ref{BoundedWeightExample} since $|\mu(k)|\leq1$ for every $k\in\mathbb{N}_0$. Hence Theorem \ref{SubsequenceResult} implies that
	$$\frac{1}{n}\sum_{k=0}^{n-1}\mu(k)T^{k^t}(x_k)-M_{\textbf{k}^t,n}^{\mu}(T)(x)=\frac{1}{n}\sum_{k=0}^{n-1}T^{k^t}(\mu(k)(x_k-x))\to0\text{ b.a.u. as }n\to\infty.$$ Since $(M_{\textbf{k}^t,n}^{\mu}(T)(x))_{n=1}^{\infty}$ converges b.a.u. as $n\to\infty$ by the above, the second conclusion follows.
	
	The above arguments show that $(M_{\textbf{k}^{t},n}^{\mu}(T)(x))_{n=1}^{\infty}$ converges to $0$ b.a.u. as $n\to\infty$ for every $x\in L_2(\mathcal{M},\tau)$, and so notably for $x\in L_1(\mathcal{M},\tau)\cap\mathcal{M}$. 
	
	Assume now that $\tau(\textbf{1})<\infty$ (noting then that $L_1(\mathcal{M},\tau)\cap\mathcal{M}=\mathcal{M}$). Then it follows by \cite[Proposition 4.2]{comli3} and the above that $\|M_{\textbf{k}^{t},n}^{\mu}(T)(x)\|_p\to0$ as $n\to\infty$ for every $x\in\mathcal{M}$ since $\|M_{\textbf{k}^{t},n}^{\mu}(T)(x)\|_\infty\leq\|x\|_\infty$ for all $n\in\mathbb{N}$. Since $\sup_{n\in\mathbb{N}}\|M_{\textbf{k}^{t},n}^{\mu}(T)\|_{L_p\to L_p}\leq1$, it follows by \cite[Lemma 2.3]{comli1} that the norm convergence occurs for every $x\in L_p(\mathcal{M},\tau)$ since $\mathcal{M}$ is dense in this space. The map $F_{T,\textbf{k}^{t}}^{\mu}:L_p(\mathcal{M},\tau)\to L_p(\mathcal{M},\tau)$ defined by $F_{T,\textbf{k}^{t}}^{\mu}(x):=\lim_{n\to\infty}M_{\textbf{k}^{t},n}^{\mu}(T)(x)$ (with the limit being with respect to $\|\cdot\|_p$) is therefore well-defined, linear, and bounded by a corollary of the Uniform Boundedness Principle. Since this (continuous) map is equal to $0$ for all elements of the dense subspace $\mathcal{M}$, it must be that $F_{T,\textbf{k}^{t}}^{\mu}(x)=0$ for all $x\in L_p(\mathcal{M},\tau)$. Since $L_p(\mathcal{M},\tau)$-norm limits and b.a.u. limits agree due to the fact that both convergences imply convergence in measure on $L_0(\mathcal{M},\tau)$ which is Hausdorff, it follows that $M_{\textbf{k}^{t},n}^{\mu}(T)(x)\to0$ b.a.u. as $n\to\infty$ for every $x\in L_p(\mathcal{M},\tau)$.
\end{proof}

%If $1<p<\infty$, then if it could be shown that $(M_{\textbf{k}^t,n}^{\mu}(T)(x))_{n=1}^{\infty}$ converges with respect to the norm of $L_p(\mathcal{M},\tau)$ for every $x\in L_p(\mathcal{M},\tau)$, then the density of $L_1(\mathcal{M},\tau)\cap\mathcal{M}$ in both $L_2(\mathcal{M},\tau)$ and $L_p(\mathcal{M},\tau)$ would imply that the limit operator is $0$ on $L_p(\mathcal{M},\tau)$ since $M_{\textbf{k}^{t},n}^{\mu}(T)(x)\to0$ b.a.u. as $n\to\infty$ for each $x\in L_1(\mathcal{M},\tau)\cap\mathcal{M}$. However, without an appropriate analogue of \cite[Theorem 1.2]{lot} for these averages, we are left without an identification of the limit operator unless $p=2$.

A similar proof to the above allows one to replace $(\mu(k))_{k=0}^{\infty}$ with a bounded invertible automatic sequence $\alpha$ by \cite[Theorem C]{ek}.  Furthermore, in that case it follows by \cite{hlrw} that one can replace $(k^t)_{k=0}^{\infty}$ by $(P(k))_{k=0}^{\infty}$, where $P(X)\in\mathbb{Z}[X]$ is a polynomial with integer coefficients such that $P(\mathbb{N}_0)\subseteq\mathbb{N}_0$.

We note that these weighted subsequential averages are a bit different than those considered in \cite{ob1}, since there the subsequence $(\alpha_{k_n})_{n=1}^{\infty}$ of $(\alpha_n)_{n=1}^{\infty}$ was used when $(k_n)_{n=0}^{\infty}$ had density $1$. However, since $\mu(k^t)=0$ for every $k\in\mathbb{N}$ when $t\geq2$, the results that one would obtain from \cite{ob1} would be trivial for the sequence $(\mu(k))_{k=0}^{\infty}$.

\subsection{Entangled Averages}

It was shown in \cite{jx} that if $S\in DS^+(\mathcal{M},\tau)$ is such that the restriction of $S$ to $L_2(\mathcal{M},\tau)$ is positive as a Hilbert-space operator (i.e. if $S\in\mathcal{B}(L_2(\mathcal{M},\tau))^{+}$), then for each $1<p<\infty$ it follows that $(S^k)_{k\in\mathbb{N}_0}$ is of strong type $(p,p)$ and $(S^k(x)-F_S(x))_{k=0}^{\infty}\in L_p(\mathcal{M},\tau;c_0)$ for every $x\in L_p(\mathcal{M},\tau)$. Applying this to Theorem \ref{GeneralBreimanErgodicTheorem} it follows that 
$$\frac{1}{n}\sum_{k=0}^{n-1}T^kS^k(x)\to F_{T}F_{S}(x)\text{ b.a.u. as }n\to\infty$$ for every $x\in L_p(\mathcal{M},\tau)$. We note that this resembles a simple form of the entangled ergodic averages considered by Liebscher in \cite{lieb} (see also \cite{ekk}). 

This can be generalized a bit more than above. To that end, we will consider averages of the form
$$\frac{1}{n}\sum_{k=0}^{n-1}T^kRS^k(x),$$ where $1<p<\infty$, $T,S\in DS^+(\mathcal{M},\tau)$ are such that there exists $d\in\mathbb{N}$ such that $S^d\in\mathcal{B}(L_2(\mathcal{M},\tau))^{+}$, $R:L_p(\mathcal{M},\tau)\to L_p(\mathcal{M},\tau)$ is a positive bounded linear operator, and $x\in L_p(\mathcal{M},\tau)$. %In the literature these are sometimes called entangled ergodic averages.

%The cases considered in this article will put a few extra restrictions on the operator $S$ that are known to hold if $1<p<\infty$ and if the restriction of $S^d$ to $L_2(\mathcal{M},\tau)$ is positive for some $d\in\mathbb{N}$. However, as far as the author is aware, % most results regarding such averages require some extra assumptions on the operators involved (at least if one is considering those in $DS^+(\mathcal{M},\tau)$), and 
%no results regarding these averages have been proven at all in the noncommutative setting.

The methods will follow the same idea as the previous sections, where we will look at a few relevant linear maps from $L_p(\mathcal{M},\tau;\ell_\infty)$ which leave $L_p(\mathcal{M},\tau;c_0)$ invariant. The main difficulty will lie in the fact that it only certain subsequences of $(S^k(x)-F_S(x))_{k=0}^{\infty}$ will be in $L_p(\mathcal{M},\tau;c_0)$ (and not necessarily the full sequence).

\begin{pro}\label{StrongTypeConstruction}
Let $(\mathcal{M},\tau)$ be a semifinite von Neumann algebra, $1\leq p<\infty$, and $T_k:L_p(\mathcal{M},\tau)\to L_p(\mathcal{M},\tau)$ be a positive bounded linear operator for every $k\in\mathbb{N}_0$. Furthermore, assume that $T=(T_k)_{k\in\mathbb{N}_0}$ is of strong type $(p,p)$ with constant $C$. Let $\overline{T}:L_p(\mathcal{M},\tau;\ell_\infty)\to L_p(\mathcal{M},\tau;\ell_\infty)$ denote the linear operator defined by
$$\overline{T}\left((x_k)_{k=0}^{\infty}\right):=(T_k(x_k))_{k=0}^{\infty}.$$ Then $\overline{T}$ is continuous with an operator norm of at most $6C$ which maps $L_p(\mathcal{M},\tau;c_0)$ to itself.
\end{pro}

\begin{proof}
Let $(x_k)_{k=0}^{\infty}\in L_p(\mathcal{M},\tau;\ell_\infty)$. By writing each $x_k$ as $x_k=y_k+iz_k$ with $y_k=\frac{x_k+x_k^*}{2}$ and $z_k=\frac{x_k-x_k^*}{2i}$, then $y_k^*=y_k$ and $z_k^*=z_k$ for every $k\in\mathbb{N}_0$, and  $$\left\|{\sup_{k\in\mathbb{N}_0}}^{+}y_k\right\|_p,\left\|{\sup_{k\in\mathbb{N}_0}}^{+}z_k\right\|_p\leq\left\|{\sup_{k\in\mathbb{N}_0}}^{+}x_k\right\|_p.$$ Hence, without loss of generality we assume that $x_k^*=x_k$ for every $k\in\mathbb{N}_0$.

Now fix $\epsilon>0$ and let $c\in L_p(\mathcal{M},\tau)^{+}$ be such that $\|c\|_p<\|{\sup_{k\in\mathbb{N}_0}}^{+}x_k\|_p+\frac{\epsilon}{3}$ and $-c\leq x_k\leq c$ for every $k\in\mathbb{N}_0$. Since each $T_k$ is positive, applying $T_k$ to each $x_k$ shows that
$$-T_k(c)\leq T_k(x_k)\leq T_k(c)\text{ for every }k\in\mathbb{N}_0.$$ Since  $c\in L_p(\mathcal{M},\tau)^{+}$ and $(T_k)_{k\in\mathbb{N}_0}$ is of strong type $(p,p)$ with constant $C$, there exists $d\in L_p(\mathcal{M},\tau)^{+}$ such that $\|d\|_p\leq C\|c\|_p$ and $0\leq T_k(c)\leq d$ for every $k\in\mathbb{N}_0.$ Therefore $$-d\leq T_k(x_k)\leq d\text{ for every }k\in\mathbb{N}_0$$ and 
$$\|d\|_p\leq C\|c\|_p< C\left(\left\|{\sup_{k\in\mathbb{N}_0}}^{+}x_k\right\|_p+\frac{\epsilon}{3(C+1)}\right)=C\left\|{\sup_{k\in\mathbb{N}_0}}^{+}x_k\right\|_p+\frac{\epsilon}{3}.$$ 
By \cite[Lemma 2.1]{lmz} it follows that $(T_k(x_k))_{k=0}^{\infty}\in L_p(\mathcal{M},\tau;\ell_\infty)$ and 
$$
\left\|{\sup_{k\in\mathbb{N}_0}}^{+}T_k(x_k)\right\|_p
\leq\|d+2d\|_p
=3\|d\|_p
< 3C\left(\left\|{\sup_{k\in\mathbb{N}_0}}^{+}x_k\right\|_p+\frac{\epsilon}{3}\right)
\leq3C\left\|{\sup_{k\in\mathbb{N}_0}}^{+}x_k\right\|_p+\epsilon.
$$
Since $\epsilon>0$ was arbitrary, it follows that $\|{\sup_{k\in\mathbb{N}_0}}^{+}T_k(x_k)\|_p\leq3C\|{\sup_{k\in\mathbb{N}_0}}^{+}x_k\|_p$.

Now let $(x_k)_{k=0}^{\infty}\in L_p(\mathcal{M},\tau;c_{00})$ and let $N\in\mathbb{N}$ be such that $x_k=0$ for $k>N$. It immediately follows that $T_k(x_k)=0$ for $k>N$ as well, meaning $(T_k(x_k))_{k=0}^{\infty}$ is a finite sequence in $L_p(\mathcal{M},\tau;\ell_\infty)$, i.e. $\overline{T}(L_p(\mathcal{M},\tau;c_{00}))\subseteq L_p(\mathcal{M},\tau;c_{00})$. Since $\overline{T}$ is continuous and since $L_p(\mathcal{M},\tau;c_{00})$ is dense in $L_p(\mathcal{M},\tau;c_0)$, it follows that
$$\overline{T}(L_p(\mathcal{M},\tau;c_0))=\overline{T}(\overline{L_p(\mathcal{M},\tau;c_{00})}\subseteq \overline{\overline{T}(L_p(\mathcal{M},\tau;c_{00}))}\subseteq L_p(\mathcal{M},\tau;c_0).$$ Therefore $(T_k(x_k))_{k=0}^{\infty}\in L_p(\mathcal{M},\tau;c_0)$ for every $(x_k)_{k=0}^{\infty}\in L_p(\mathcal{M},\tau;c_0)$.
\end{proof}

The following is a mild improvement of \cite[Proposition 7.3]{hjx}. In particular, using the same methods as above we can reduce the operator norm of the operator considered from $8\|R\|_{L_p\to L_p}$ to $6\|R\|_{L_p\to L_p}$, show that the map takes $c_0$-valued spaces to $c_0$-valued spaces, and prove that it remains true when $R$ maps between different von Neumann algebras. Although we state it for semifinite von Neumann algebras, this would also extend to the general von Neumann algebra case for the same reasons described in Remark \ref{GeneralvonNeumannAlgebras} (which is what was considered in \cite{hjx}).

\begin{pro}\label{SingleOperatorMaximalInequality}
Let $(\mathcal{M},\tau)$ and $(\mathcal{N},\rho)$ be semifinite von Neumann algebras, fix $1\leq p<\infty$, and let $R:L_p(\mathcal{M},\tau)\to L_p(\mathcal{N},\rho)$ be a positive bounded linear operator. Define $\overline{R}:L_p(\mathcal{M},\tau;\ell_\infty)\to L_p(\mathcal{N},\rho;\ell_\infty)$ by $$\overline{R}((x_k)_{k=0}^{\infty}):=(R(x_k))_{k=0}^{\infty}.$$ Then $\overline{R}$ is a bounded linear operator with a norm of at most $6\|R\|_{L_p\to L_p}$. Furthermore, $\overline{R}$ maps $L_p(\mathcal{M},\tau;c_0)$ to $L_p(\mathcal{N},\rho;c_0)$.
\end{pro}
\begin{proof}
Let $(x_k)_{k=0}^{\infty}\in L_p(\mathcal{M},\tau;\ell_\infty)$ and fix $\epsilon>0$. We first assume that $x_k^*=x_k$ for every $k\in\mathbb{N}_0$. By \cite[Lemma 2.1(i)]{lmz} we find that $(x_k)_{k=0}^{\infty}\in L_p(\mathcal{M},\tau;\ell_\infty)$ being a self-adjoint sequence implies that there exists $c\in L_p(\mathcal{M},\tau)^{+}$ with $$\|c\|_p<\left\|{\sup_{k\in\mathbb{N}_0}}^{+}x_k\right\|_p+\frac{\epsilon}{3(\|R\|_{L_p\to L_p}+1)}\text{ and }-c\leq x_k\leq c\text{ for every }k\in\mathbb{N}_0.$$ 

Since $R$ is a positive operator, it follows that $R(c)\in L_p(\mathcal{N},\rho)^{+}$ with
$$-R(c)\leq R(x_k)\leq R(c)\text{ for every }k\in\mathbb{N}_0$$ and, since $R$ is also bounded, that
\begin{align*}
	\|R(c)\|_p
	\leq&\|R\|_{L_p\to L_p}\|c\| \\
	<&\|R\|_{L_p\to L_p}\left(\left\|{\sup_{k\in\mathbb{N}_0}}^{+}x_k\right\|_p+\frac{\epsilon}{3(\|R\|_{L_p\to L_p}+1)}\right) \\
	\leq&\|R\|_{L_p\to L_p}\left\|{\sup_{k\in\mathbb{N}_0}}^{+}x_k\right\|_p+\frac{\epsilon}{3}.
\end{align*}
By \cite[Lemma 2.1(ii)]{lmz} it follows that $(R(x_k))_{k=0}^{\infty}\in L_p(\mathcal{N},\rho;\ell_\infty)$ with $$\left\|{\sup_{k\in\mathbb{N}_0}}^{+}R(x_k)\right\|_p\leq 3\|R\|_{L_p\to L_p}\left\|{\sup_{k\in\mathbb{N}_0}}^{+}x_k\right\|_p+\epsilon,$$ which since $\epsilon>0$ was arbitrary it follows that the same holds with $\epsilon=0$. 

For general $(x_k)_{k=0}^{\infty}\in L_p(\mathcal{M},\tau;\ell_\infty)$ write $x_k=y_k+iz_k$ with $y_k=\frac{x_k+x_k^*}{2}$ and $z_k=\frac{x_k-x_k^*}{2i}$ for every $k\in\mathbb{N}_0$. Applying the above to $(y_k)_{k=0}^{\infty}$ and $(z_k)_{k=0}^{\infty}$, using the fact that the $L_p(\mathcal{M},\tau;\ell_\infty)$-norms of $(y_k)_{k=0}^{\infty}$ and $(z_k)_{k=0}^{\infty}$ are smaller than that of $(x_k)_{k=0}^{\infty}$, and using the triangle inequality of the norm of $(x_k)_{k=0}^{\infty}$ one finds that $\left\|{\sup_{k\in\mathbb{N}_0}}^{+}R(x_k)\right\|_p\leq6\|R\|_{L_p\to L_p}\left\|{\sup_{k\in\mathbb{N}_0}}^{+}x_k\right\|_p$ for every $(x_k)_{k=0}^{\infty}\in L_p(\mathcal{M},\tau;\ell_\infty)$.

A similar proof to that of Proposition \ref{StrongTypeConstruction} shows that $L_p(\mathcal{M},\tau;c_0)$ is mapped to $L_p(\mathcal{N},\rho;c_0)$ under $\overline{R}$.
%Since $(x_k)_{k=0}^{\infty}\in L_p(\mathcal{M},\tau;\ell_\infty)$ with $x_k^*=x_k$ for every $k\in\mathbb{N}_0$ was arbitrary, it follows that $(R)_{k\in\mathbb{N}_0}$ is of strong type $(p,p)$ with constant $6\|R\|_{L_p\to L_p}$ as in the beginning of the proof of Proposition \ref{StrongTypeConstruction}. 
\end{proof}

One application of this is the following.

\begin{cor}
Let $(\mathcal {M},\tau)$ be a semifinite von Neumann algebra, fix $1\leq p<\infty$, and suppose that $\mathcal{N}\subset\mathcal{M}$ is a von Neumann subalgebra such that $\tau|_{\mathcal{N}}$ is also semifinite. Let $\mathbb{E}:\mathcal{M}\to\mathcal{N}$ denote the corresponding conditional expectation (as well as its bounded extension from $L_p(\mathcal{M},\tau)$ to $L_p(\mathcal{N},\tau|_{\mathcal{N}})$). Then if $(x_k)_{k=0}^{\infty}\in L_p(\mathcal{M},\tau;\ell_\infty)$ it follows that
$$(\mathbb{E}(x_k))_{k=0}^{\infty}\in L_p(\mathcal{N},\tau|_{\mathcal{N}};\ell_\infty).$$ Furthermore, the same conclusion holds when each instance of $\ell_\infty$ is replaced by $c_0$.
\end{cor}

We will now prove the following entangled ergodic theorem.

\begin{teo}\label{EntangledAverages2}
Let $(\mathcal{M},\tau)$, $(\mathcal{N},\rho)$ be semifinite von Neumann algebras and fix $1<p<\infty$. Let $T\in DS^+(\mathcal{N},\rho)$, $R:L_p(\mathcal{M},\tau)\to L_p(\mathcal{N},\rho)$ be a positive bounded linear operator, and let $S\in DS^+(\mathcal{M},\tau)$ be such that the restriction of $S^d$ to $L_2(\mathcal{M},\tau)$ is positive as a Hilbert space operator for some $d\in\mathbb{N}$ (i.e. $S^d\in\mathcal{B}(L_2(\mathcal{M},\tau))^{+}$).	

Then for every $x\in L_p(\mathcal{M},\tau)$ it follows that the entangled averages of $T$, $R$, and $S$ satisfy
$$\left(\frac{1}{n}\sum_{k=0}^{n-1}T^kRS^k(x)-\frac{1}{d}\sum_{j=0}^{d-1}F_{T^d}T^jRF_{S^d}S^j(x)\right)_{n=1}^{\infty}\in L_p(\mathcal{N},\rho;c_0).$$  	
Furthermore, if $(x_k)_{k=0}^{\infty}\in L_p(\mathcal{M},\tau;\ell_\infty)$ and $x\in L_p(\mathcal{M},\tau)$ are such that $(x_k-x)_{k=0}^{\infty}\in L_p(\mathcal{M},\tau;c_0)$, then
$$\left(\frac{1}{n}\sum_{k=0}^{n-1}T^kRS^k(x_k)-\frac{1}{d}\sum_{j=0}^{d-1}F_{T^d}T^jRF_{S^d}S^j(x)\right)_{n=1}^{\infty}\in L_p(\mathcal{N},\rho;c_0).$$
Consequently, $$\frac{1}{n}\sum_{k=0}^{n-1}T^kRS^k(x), \ \frac{1}{n}\sum_{k=0}^{n-1}T^kRS^k(x_k)\to\frac{1}{d}\sum_{j=0}^{d-1}F_{T^d}T^{j}RF_{S^d}S^{j}(x)$$ b.a.u. as $n\to\infty$ for every $x\in L_p(\mathcal{M},\tau)$ and $(x_k)_{k=0}^{\infty}\in L_p(\mathcal{M},\tau;\ell_\infty)$ such that $(x_k-x)_{k=0}^{\infty}\in L_p(\mathcal{M},\tau;c_0)$.
\end{teo}

\begin{proof}
We will first prove the claim for the standard averages, and will deal with the case of replacing $T^kRS^k(x)$ with $T^kRS^k(x_k)$ for every $k\in\mathbb{N}_0$ afterwards. As such, fix $x\in L_p(\mathcal{M},\tau)$. 

Write $a\mathbb{N}_0+b:=\{an+b:n\in\mathbb{N}_0\}$ for every $a,b\in\mathbb{N}_0$ and $\chi_E$ for the characteristic function of a set $E\subseteq\mathbb{N}_0$. One can find that
\begin{align}\label{EqualityI}
\left(S^k(x)-\sum_{j=0}^{d-1}\chi_{d\mathbb{N}_0+j}(k)F_{S^d}(S^j(x))\right)_{k=0}^{\infty}
=&\sum_{j=0}^{d-1}\left(\chi_{d\mathbb{N}_0+j}(k)\left(S^k(x)-F_{S^d}(S^j(x))\right)\right)_{k=0}^{\infty}
\end{align}
as elements of $L_p(\mathcal{M},\tau;\ell_\infty)$. Since $S^d$ is positive as an operator on the Hilbert space $L_2(\mathcal{M},\tau)$, it follows that $$(S^{dk+j}(x)-F_{S^d}(S^j(x)))_{k=0}^{\infty}=((S^d)^{k}(S^j(x))-F_{S^d}(S^j(x)))_{k=0}^{\infty}\in L_p(\mathcal{M},\tau;c_0)$$ for every $j\in\{0,1,...,d-1\}$ by \cite[Theorem 6.7]{jx}. From this on sees that every term on the right hand side of Equation (\ref{EqualityI}) above is in $L_p(\mathcal{M},\tau;c_0)$ as well by placing $0$'s between appropriate terms in the factorization as elements of $L_p(\mathcal{M},\tau;\ell_\infty)$. Hence, the left hand side of (\ref{EqualityI}) is in $L_p(\mathcal{M},\tau;c_0)$.

By applying Proposition \ref{SingleOperatorMaximalInequality} followed by Theorem \ref{GeneralBreimanErgodicTheorem} and then applying the linearity of each operator it follows that
\begin{align}\label{EqualityII}
\left(\frac{1}{n}\sum_{k=0}^{n-1}T^kRS^k(x)-\sum_{j=0}^{d-1}\frac{1}{n}\sum_{k=0}^{n-1}\chi_{d\mathbb{N}_0+j}(k)T^k(RF_{S^d}S^j(x))\right)_{n=1}^{\infty}\in L_p(\mathcal{N},\rho;c_0).
\end{align}

Fix $j\in\{0,1,...,d-1\}$ and let $m_j(n)=\sup\{m\in\mathbb{N}_0:dm+j\leq n-1\}$ for every $n\in\mathbb{N}$. Note that $m_j(n)\to\infty$ and $\frac{m_j(n)}{n}\to\frac{1}{d}$ both as $n\to\infty$. Then for every $n\in\mathbb{N}$ we find
\begin{align*}
\frac{1}{n}\sum_{k=0}^{n-1}\chi_{d\mathbb{N}_0+j}(k)T^k(RF_{S^d}S^j(x))
=&\frac{m_j(n)+1}{n}\frac{1}{m_j(n)+1}\sum_{m=0}^{m_j(n)}(T^d)^{m}(T^jRF_{S^d}S^j(x)) \\
=&\frac{m_j(n)+1}{n}M_{m_j(n)+1}(T^d)(T^jRF_{S^d}S^j(x)).
\end{align*}

Hence it follows by Corollary \ref{NumberLimits} that
$$\left(\frac{1}{n}\sum_{k=0}^{n-1}\chi_{d\mathbb{N}_0+j}(k)T^k(RF_{S^d}S^j(x))-\frac{1}{d}F_{T^d}T^jRF_{S^d}S^j(x)\right)_{n=1}^{\infty}\in L_p(\mathcal{N},\rho;c_0).$$
Using this with (\ref{EqualityII}) implies that
$$\left(\frac{1}{n}\sum_{k=0}^{n-1}T^kRS^k(x)-\frac{1}{d}\sum_{j=0}^{d-1}F_{T^d}T^jRF_{S^d}S^j(x)\right)_{n=1}^{\infty}\in L_p(\mathcal{N},\rho;c_0),$$ 
finishing the proof of the first claim.

If we can show that $S=(S^k)_{k\in\mathbb{N}_0}$ is of strong type $(p,p)$, then applying Proposition \ref{StrongTypeConstruction} to $(x_k-x)_{k=0}^{\infty}\in L_p(\mathcal{M},\tau;c_0)$ followed by Proposition \ref{SingleOperatorMaximalInequality}, Theorem \ref{GeneralBreimanErgodicTheorem}, and the first part of this proof would prove the second claim.

We will now show that $(S^k)_{k\in\mathbb{N}_0}$ is of strong type $(p,p)$. Fix $y\in L_p(\mathcal{M},\tau)^{h}$. Since $S^d$ is a positive as a Hilbert space operator on $L_2(\mathcal{M},\tau)$ it is known that $((S^{md}))_{m\in\mathbb{N}_0}=((S^d)^{m})_{m\in\mathbb{N}_0}$ is of strong type $(p,p)$ with constant $C_p$. Hence, for each $j\in\{0,1,...,d-1\}$ there exists $c_j\in L_p(\mathcal{M},\tau)^{+}$ such that $$\|c_j\|_p\leq C_p\|S^j(y)\|_p\text{ and }-c_j\leq S^{md}(S^j(y))\leq c_j\text{ for every }m\in\mathbb{N}_0.$$
Since $c_j\leq\sum_{\ell=0}^{d-1}c_\ell$ for every $j\in\mathbb{N}_0$ it follows that $$-\sum_{\ell=0}^{d-1}c_{\ell}\leq (S^d)^{m}(S^j(y))=S^{md+j}(y)\leq\sum_{\ell=0}^{d-1}c_\ell\text{ for every }m\in\mathbb{N}_0\text{ and }j\in\{0,...,d-1\},$$ and 
$$\left\|\sum_{\ell=0}^{d-1}c_\ell\right\|_p\leq\sum_{\ell=0}^{d-1}\|c_\ell\|_p\leq\sum_{\ell=0}^{d-1}C_p\|S^\ell(y)\|_p\leq\sum_{\ell=0}^{d-1}C_p\|y\|_p=dC_p\|y\|_p.$$ Since every $k\in\mathbb{N}_0$ can be written as $k=md+j$ for some $m\in\mathbb{N}_0$ and $j\in\{0,...,d-1\}$, it follows that $(S^k(y))_{k\in\mathbb{N}_0}\in L_p(\mathcal{M},\tau;\ell_\infty)$ with a norm of at most $3dC_p\|y\|_p$ by \cite[Lemma 2.1]{lmz}.

If $w\in L_p(\mathcal{M},\tau)$, then write $w=y+iz$ for some $y,z\in L_p(\mathcal{M},\tau)^{h}$ with $\|y\|_p,\|z\|_p\leq\|w\|_p$. Then
$$
\left\|{\sup_{k\in\mathbb{N}_0}}^{+}S^k(w)\right\|_p
\leq\left\|{\sup_{k\in\mathbb{N}_0}}^{+}S^k(y)\right\|_p+\left\|{\sup_{k\in\mathbb{N}_0}}^{+}S^k(z)\right\|_p 
\leq3dC_p\|y\|_p+3dC_p\|z\|_p 
%\leq3dC_p\|x\|_p+3dC_p\|x\|_p
\leq6dC_p\|w\|_p.
$$
Hence $(S^k(w))_{k=0}^{\infty}\in L_p(\mathcal{M},\tau;\ell_\infty)$ with a norm of at most $6dC_p\|w\|_p$, which implies that $(S^k)_{k\in\mathbb{N}_0}$ is of strong type $(p,p)$ with constant $6dC_p$.
\end{proof}

An important case of the operator $S$ above is when the restriction of $S$ to $L_2(\mathcal{M},\tau)$ is self-adjoint as a Hilbert space operator. In this case, the result above applies by using $d=2$.

One example of this type of operator can be defined as follows. If $\Phi:\mathcal{M}\to\mathcal{M}$ is a normal $\tau$-preserving $*$-automorphism (writing $\Phi$ for its isometric extensions to $L_p(\mathcal{M},\tau)$ for $1\leq p\leq\infty$ as well), and if $\nu:P(\mathbb{Z})\to[0,1]$ is a probability measure on the power set $P(\mathbb{Z})$ of the integers $\mathbb{Z}$ with $\nu(\{-k\})=\nu(\{k\})$ for every $k\in\mathbb{Z}$, then the operator $S:L_1(\mathcal{M},\tau)+\mathcal{M}\to L_1(\mathcal{M},\tau)+\mathcal{M}$ defined by
$$S(x):=\sum_{k\in\mathbb{Z}}\nu(\{k\})\Phi^k(x)\text{ for }x\in L_1(\mathcal{M},\tau)+\mathcal{M}$$ is an element of $DS^+(\mathcal{M},\tau)$ which is self-adjoint when restricted to $L_2(\mathcal{M},\tau)$.

The assumption of $S^d$ being positive is not the only assumption on $S$ for which a version result would apply. In general, what was needed was that there exists $d\in\mathbb{N}$ such that $(S^{kd}(x)-F_{S^d}(x))_{k=0}^{\infty}\in L_p(\mathcal{M},\tau;c_0)$ for every $x\in L_p(\mathcal{M},\tau)$ and $(S^{kd})_{k\in\mathbb{N}_0}$ is of strong type $(p,p)$. See \cite{be} for some other assumptions on $S$ which satisfy this assumption.

An explicit example of an operator satisfying the assumption of \cite{be} is as follows. Let $(\mathcal{M},\tau)$ be a semifinite von Neumann algebra, let $\mathcal{M}_1,...,\mathcal{M}_m$ be von Neumann subalgebras of $\mathcal{M}$ such that $\tau|_{\mathcal{M}_j}$ is semifinite on $\mathcal{M}_j$ for every $j\in\{1,...,m\}$, and let $\mathbb{E}_j:\mathcal{M}\to\mathcal{M}_j$ denote the corresponding conditional expectation. Then it was shown in \cite{be} that the operator $S\in DS^+(\mathcal{M},\tau)$ defined by $$S(x):=\mathbb{E}_1\cdots\mathbb{E}_m(x)\text{ for }x\in L_1(\mathcal{M},\tau)+\mathcal{M}$$ satisfies the condition above with $d=1$ with $F_S$ being the conditional expectation of $\mathcal{M}$ onto the subalgebra $\bigcap_{j=1}^{m}\mathcal{M}_j$.

\begin{comment}
If one only had that $(x_k)_{k=0}^{\infty}\subset L_p(\mathcal{M},\tau)$ is a sequence such that $x_k\to 0$ b.a.u. as $k\to\infty$ (without necessarily knowing if $(x_k)_{k=0}^{\infty}\in L_p(\mathcal{M},\tau;c_0)$), then it is unknown if this result is still true unless $T:\mathcal{M}\to\mathcal{M}$ is a normal $\tau$-preserving $*$-homorphism.

The main issue lies in the fact that, if $e\in\mathcal{P(M)}$ is such that $\|ex_{k}e\|_\infty\to0$ as $k\to\infty$, then in general there is nothing guaranteeing that there is a projection $f\in\mathcal{P(M)}$ with $\|fT(x_k)f\|_\infty\to0$. 

However, if $T$ is a $*$-homomorphism described above, then $T(e)\in\mathcal{P(M)}$ satisfies $\tau(T(e)^{\perp})=\tau(e^{\perp})$ and $$\|T(e)T(x_k)T(e)\|_\infty=\|T(e(x_k)e)\|_\infty=\|e(x_k)e\|_\infty\to0$$ as $k\to\infty$.
\end{comment}

\subsection{$T$-admissible processes}

We end this article with a discussion of one more method of finding sequences in $L_p(\mathcal{M},\tau;c_0)$ and some consequences of these relating to $T$-admissible sequences as studied in \cite{comli1}, \cite{comli2}, and \cite{ob5}.

Following \cite{comli2}, call a sequence $(x_k)_{k=0}^{\infty}\subset L_p(\mathcal{M},\tau)$ strongly $p$-bounded if $$\sup_{k\in\mathbb{N}_0}\|x_k\|_p<\infty.$$

It was shown in \cite[Remark 3.1(2)]{dj} that if $(x_k)_{k=0}^{\infty}\in L_p(\mathcal{M},\tau;\ell_\infty)$, then
	$$\sup_{k\in\mathbb{N}_0}\|x_k\|_p\leq\left\|{\sup_{k\in\mathbb{N}_0}}^{+}x_k\right\|_p.$$	
Hence, every element of $L_p(\mathcal{M},\tau;\ell_\infty)$ is strongly $p$-bounded.

\begin{pro}\label{IncreasingSequences}
	Let $(\mathcal{M},\tau)$ be a semifinite von Neumann algebra, $1\leq p<\infty$, and $(x_k)_{k=0}^{\infty}\subset L_p(\mathcal{M},\tau)^{h}$ be a sequence which is strongly $p$-bounded and is increasing with respect to the ordering $\leq$ on $L_0(\mathcal{M},\tau)^{h}$. Then there exists $x\in L_p(\mathcal{M},\tau)^h$ such that
	$$(x_k-x)_{k=0}^{\infty}\in L_p(\mathcal{M},\tau;c_0).$$
\end{pro}
\begin{proof}
	Since $x_k\leq x_{k+1}$ for every $k\in\mathbb{N}_0$ and $\sup_{k\in\mathbb{N}_0}\|x_k\|_p<\infty$ both by assumption, we find that $0\leq x_k-x_0\leq x_{k+1}-x_0$ for every $k\in\mathbb{N}_0$ and $\sup_{k\in\mathbb{N}_0}\|x_k-x_0\|_p\leq\|x_0\|_p+\sup_{k\in\mathbb{N}_0}\|x_k\|_p<\infty$. Therefore it follows by \cite[Proposition 3.2 and Theorem 5.15]{ddp} that there exists $y\in L_p(\mathcal{M},\tau)^+$ such that $0\leq x_k-x_0\leq y$ for every $k\in\mathbb{N}_0$ and $\|y-(x_k-x_0)\|_p\to0$ as $k\to\infty$; in particular, $y=\sup_{k\in\mathbb{N}_0}(x_k-x_0)$ (with respect to the ordering $\leq$ on $L_0(\mathcal{M},\tau)^+$). If $x:=y+x_0$ then $x\geq x_k$ (so that $x-x_k\geq0$) for every $k\in\mathbb{N}_0$ and $$\|x-x_k\|_p=\|(y+x_0)-x_k\|_p=\|y-(x_k-x_0)\|_p\to0\text{ as }k\to\infty.$$
	
	Fix $\epsilon>0$.
	
	It follows from $\|x-x_k\|_p\to0$ that there exists $k^{\prime}\in\mathbb{N}_0$ such that $k\geq k^{\prime}$ implies that $\|x-x_{k}\|_p\leq\epsilon$. Since $(x_k)_{k=0}^{\infty}$ is increasing with respect to the ordering on $L_0(\mathcal{M},\tau)^h$ it follows that $(x-x_k)_{k=0}^{\infty}$ is decreasing in $L_0(\mathcal{M},\tau)^+$ since $x_k\leq x$ for all $k$. Therefore 
	$$0\leq x-x_k\leq x-x_{k^{\prime}}\text{ for every }k\geq k^{\prime}$$ and $\max\{\|0\|_p,\|x-x_{k^{\prime}}\|_p\}\leq\epsilon$. Since $\epsilon>0$ was arbitrary, the conclusion follows from \cite[Corollary 2.3]{lmz}.
\end{proof}

Fix $1\leq p<\infty$ and let $T\in DS^+(\mathcal{M},\tau)$ be  induced by a normal $\tau$-preserving $*$-automorphism of $\mathcal{M}$ (so that $T^{-1}\in DS^+(\mathcal{M},\tau)$ as well). A sequence $(x_k)_{k=0}^{\infty}\subset L_p(\mathcal{M},\tau)^h$ is called \textit{$T$-admissible} if
$$T(x_k)\leq x_{k+1}\text{ for every }k\in\mathbb{N}_0.$$

If $(x_k)_{k=0}^{\infty}$ is a strongly $p$-bounded $T$-admissible process, then it follows from \cite[Lemma 4.1]{ob5} that there exists an operator $w\in L_p(\mathcal{M},\tau)^h$ such that $x_k\leq T^k(w)$ for every $k\in\mathbb{N}_0$ and that $\|x_k-T^{k}(w)\|_p\to0$ as $k\to\infty$. This is then used to prove ergodic theorems about $T$-admissible processes similar to those considered earlier in this article. However, the proof that $w$ exists is the same as for why $x$ exists in the previous result, as we will now see.

\begin{cor}\label{TAdmissible}
Let $(\mathcal{M},\tau)$ be a semifinite von Neumann algebra, $1\leq p<\infty$, and $T:L_p(\mathcal{M},\tau)\to L_p(\mathcal{M},\tau)$ be induced by a normal $\tau$-preserving $*$-automorphism of $\mathcal{M}$. If $(x_k)_{k=0}^{\infty}\subset L_p(\mathcal{M},\tau)^h$ is a strongly $p$-bounded $T$-admissible process, then there exists $w\in L_p(\mathcal{M},\tau)^h$ such that
$$(T^{-k}(x_k)-w)_{k=0}^{\infty}\in L_p(\mathcal{M},\tau;c_0).$$
Furthermore, it follows that
$$\left(\frac{1}{n}\sum_{k=0}^{n-1}x_k-F_T(w)\right)_{n=1}^{\infty}\in L_p(\mathcal{M},\tau;c_0),$$ and so $\frac{1}{n}\sum_{k=0}^{n-1}x_k\to F_T(w)$ b.a.u. as $n\to\infty$.	
\end{cor}
\begin{proof}
If $k\in\mathbb{N}_0$, then the positivity of $T^{-(k+1)}$ and $T(x_k)\leq x_{k+1}$ implies that
$$T^{-(k+1)}(x_{k+1})\geq T^{-(k+1)}(T(x_k))=T^{-k-1+1}(x_k)=T^{-k}(x_k).$$
Therefore $(T^{-k}(x_k))_{k=0}^{\infty}$ is an increasing sequence in $L_p(\mathcal{M},\tau)^h$ with $$\sup_{k\in\mathbb{N}_0}\|T^{-k}(x_k)\|_p=\sup_{k\in\mathbb{N}_0}\|x_k\|_p<\infty,$$ from which the first conclusion follows from Proposition \ref{IncreasingSequences}.

From this, since
$$\frac{1}{n}\sum_{k=0}^{n-1}T^k(T^{-k}(x_k)-w)=\frac{1}{n}\sum_{k=0}^{n-1}x_k-M_n(T)(w) \text{ for every }n\in\mathbb{N},$$ and since $(M_n(T)(w)-F_T(w))_{n=1}^{\infty}\in L_p(\mathcal{M},\tau;c_0)$ by \cite[Theorem 6.3]{jx} the second conclusion follows by Theorem \ref{GeneralBreimanErgodicTheorem}.
\end{proof}

%\begin{rem}
The b.a.u. convergence of weighted averages of $T$-admissible processes of operators from the noncommutative $L_p$-spaces associated to a semifinite von Neumann algebras was studied in \cite{ob5}. Indeed, \cite[Theorem 4.2]{ob5} states that, if $\alpha=(\alpha_k)_{k=0}^{\infty}\in W_\infty$ and $(M_n^{\alpha}(T)(y))_{n=1}^{\infty}$ converges b.a.u. as $n\to\infty$ for every $y\in L_p(\mathcal{M},\tau)$, then for any increasing strongly $p$-bounded sequence $(x_k)_{k=0}^{\infty}$ in $L_p(\mathcal{M},\tau)$ the averages $(\frac{1}{n}\sum_{k=0}^{n-1}\alpha_k T^k(x_k))_{n=1}^{\infty}$ converge b.a.u. as $n\to\infty$ to the same limit as $(\frac{1}{n}\sum_{k=0}^{n-1}\alpha_k T^k(x))_{n=1}^{\infty}$ when $x=\sup_{k\in\mathbb{N}_0}x_k$, which is exactly the result of applying Proposition \ref{IncreasingSequences} to Theorems \ref{WeightedDirect} and \ref{WeightedGeneralized}. However, the results in \cite{ob5} are stronger than those of this article in the following ways.

First off, if $1<p,q<\infty$ with $\frac{1}{p}+\frac{1}{q}=1$ then one can also consider the weight $\alpha$ to be in $W_q$, compared to the assumption made in Theorem \ref{WeightedGeneralized} above which requires them to be in the smaller space $W_s$ for some $s>q$. 

If $p=1$ then bounded weights are still needed. However, \cite[Theorem 4.2]{ob5} shows that one can consider any operator $T\in DS^+(\mathcal{M},\tau)$ instead of requiring $T$ to be induced by an ergodic $*$-homomorphism of $\mathcal{M}$ as in Theorems \ref{DirectBreimanErgodicTheorem} and \ref{WeightedDirect}. Furthermore, one can also take $(\mathcal{M},\tau)$ to be any semifinite von Neumann algebra instead of requiring $\tau(\textbf{1})<\infty$.

Next, \cite{ob5} also shows that Wiener-Wintner type results can be obtained in both of these cases as well, compared to the one sequence at a time approach taken by Theorems \ref{WeightedDirect} and \ref{WeightedGeneralized} in the present article. A Wiener-Wintner type result holds for $T$ acting on $L_p(\mathcal{M},\tau)$ with respect to the set $\mathcal{W}\subset W_q$ if for any $y\in L_p(\mathcal{M},\tau)$ and $\epsilon>0$ one can find a projection $e\in\mathcal{P(M)}$ with $\tau(e^\perp)\leq\epsilon$ such that $(eM_n^{\alpha}(T)(y)e)_{n=1}^{\infty}$ converges in $\mathcal{M}$ as $n\to\infty$ for every $\alpha\in\mathcal{W}$. Notably, the projection $e$ is independent of any particular weight $\alpha\in\mathcal{W}$, depending only possibly on $\mathcal{W}$ as a whole. 
The results in \cite[Section 4]{ob5} show that, in the case where this holds, then the same can be said if one replaces each $T^k(x)$ in these averages by $T^k(x_k)$ when $(x_k)_{k=0}^{\infty}$ is a strongly $p$-bounded increasing sequence in $L_p(\mathcal{M},\tau)^{h}$. Furthermore, the limit of these averages agrees with that of $(M_n^{\alpha}(T)(x))_{n=1}^{\infty}$ for every $\alpha\in\mathcal{W}$, where $x=\sup_{k\in\mathbb{N}_0}x_k$.

The results in \cite{ob5} only focused on studying weighted averages, while the results of this article also apply to subsequential and entangled averages as well. The assumptions made on $(x_k)_{k=0}^{\infty}$ in this article are more general than those considered in \cite{ob5}, as is shown by Proposition \ref{IncreasingSequences} and Corollary \ref{TAdmissible}.

In the situation where the results of both \cite{ob5} and the present article would apply, the main advantage of the setup considered in this paper is that b.a.u. convergence of the weighted averages occurs as a consequence of the sequence of averages being in $L_p(\mathcal{M},\tau;c_0)$ when $1<p<\infty$ (which, as previously stated, is a stronger type of convergence). %As such, the methods of this article can be generalized to many other situations (as was discussed throughout the article).
%\end{rem}

\end{document}